\documentclass[12pt]{amsart}

\usepackage{latexsym}
\usepackage{amsfonts}
\usepackage{amssymb}
\usepackage{amsmath}
\usepackage{mathrsfs}
\input xypic
\usepackage{tikz}
\usepackage{hyperref}

\newcommand{\be}{\begin{equation}}
\newcommand{\ee}{\end{equation}}
\newcommand{\bea}{\begin{eqnarray}}
\newcommand{\eea}{\end{eqnarray}}
\newcommand{\bean}{\begin{eqnarray*}}
\newcommand{\eean}{\end{eqnarray*}}
\newcommand{\brray}{\begin{array}}
\newcommand{\erray}{\end{array}}

\newtheorem{dfn}{Definition}[section]
\newtheorem{thm}[dfn]{Theorem}
\newtheorem{lmma}[dfn]{Lemma}
\newtheorem{ppsn}[dfn]{Proposition}
\newtheorem{crlre}[dfn]{Corollary}
\newtheorem{xmpl}[dfn]{Example}
\newtheorem{rmrk}[dfn]{Remark}

\newcommand{\bdfn}{\begin{dfn}\rm}
\newcommand{\bthm}{\begin{thm}}
\newcommand{\blmma}{\begin{lmma}}
\newcommand{\bppsn}{\begin{ppsn}}
\newcommand{\bcrlre}{\begin{crlre}}
\newcommand{\bxmpl}{\begin{xmpl}}
\newcommand{\brmrk}{\begin{rmrk}\rm}

\newcommand{\edfn}{\end{dfn}}
\newcommand{\ethm}{\end{thm}}
\newcommand{\elmma}{\end{lmma}}
\newcommand{\eppsn}{\end{ppsn}}
\newcommand{\ecrlre}{\end{crlre}}
\newcommand{\exmpl}{\end{xmpl}}
\newcommand{\ermrk}{\end{rmrk}}

\author{S. Sundar\\ }
\title{Cuntz-Li relations, Inverse semigroups and groupoids }

\email{sundarsobers@gmail.com}
\address{Indian Statistical Institute, Delhi.}
\begin{document}
\begin{abstract}
  In this paper we show that the universal $C^{*}$-algebra satisfying the Cuntz-Li relations is generated by an inverse semigroup of partial isometries. We apply Exel's theory of tight representations to this inverse semigroup. We identify the universal $C^{*}$-algebra as the $C^{*}$-algebra of the tight groupoid associated to the inverse semigroup.
\end{abstract}
\maketitle
\section{Introduction}
 Let $R$ be an integral domain with only finite quotients. Assume that $R$ is not a field and let $K$ be its field of fractions. We denote the set of non-zero elements in $R$ (resp. $K$) by $R^{\times}$ (resp. $K^{\times}$). In \cite{Cuntz-Li}, Cuntz and Li studied the $C^{*}$-algebra, denoted $\mathfrak{A}_{r}[R]$, on $\ell^{2}(R)$ generated by the isometries induced by the multiplication and addition operations of the ring $R$. They showed that it is simple and purely infinite. It was also shown that this $C^{*}$-algebra is the universal $C^{*}$-algebra generated by isometries satisfying the relations reflecting the semigroup multiplication in $R \rtimes R^{\times}$  and one more important relation satisfied by the range projections. Also it was shown that $\mathfrak{A}_{r}[R]$ is Morita-equivalent to a crossed product of the form $C_{0}(\mathcal{R}) \rtimes (K \rtimes K^{\times})$ where $\mathcal{R}$ is a locally compact Hausdorff space. For $R=\mathbb{Z}$, $\mathcal{R}=\mathbb{A}_{f}$ is  the space of finite adeles. Alternate approaches to the algebra $\mathfrak{A}_{r}[R]$ were considered in \cite{Quigg-Landstad}, \cite{Ex2}, and \cite{Sundar1}. 

 In \cite{Quigg-Landstad}, the situation in \cite{Cuntz-Li} was abstracted. Consider a semidirect product $N \rtimes H$ and a normal subgroup $M$ of $N$. Let $P:=\{a \in H:aMa^{-1} \subset M \}$. Then $P$ is a semigroup. In \cite{Quigg-Landstad},under certain hypotheses regarding the pair $(G=N \rtimes H,M)$, the crossed product algebra $C_{0}(\overline{N})\rtimes G$ was considered. Here $\overline{N}$ is the profinite completion of $N$ with respect to the group topology induced by the neighbourhood base $\{aMa^{-1}\}_{a\in H}$ at the identity. Let $\overline{M}$ be the closure of $M$ in $\overline{N}$. In \cite{Quigg-Landstad}, it was shown that the crossed product algebra $C_{0}(\overline{N})\rtimes G$ is Morita-equivalent to the  $C^{*}$-algebra of the groupoid $\overline{N} \rtimes G|_{\overline{M}}$.  In \cite{Quigg-Landstad}, It was shown that when $H$ is abelian, $C^{*}(\overline{N} \rtimes G|_{\overline{M}})$ is the universal $C^{*}$-algebra generated by isometries satisfying the relations reflecting the semigroup multiplication in $M \rtimes P$ and one more important relation among the range projections.  They also obtained sufficient conditions which will ensure that the reduced $C^{*}$-algebra $C^{*}_{red}(\overline{N} \rtimes G|_{\overline{M}})$ is simple and purely infinite.

Our objective in this paper is to weaken the hypothesis that $H$ is abelian. Instead we assume  $H=PP^{-1}=P^{-1}P$. This allows us to  consider pairs like $(\mathbb{Q}^{n} \rtimes GL_{n}(\mathbb{Q}),\mathbb{Z}^{n})$. Also we start with the universal $C^{*}$-algebra, denoted $\mathfrak{A}[N \rtimes H,M]$, generated by isometries satisfying the Cuntz-Li relations (See Defn. \ref{Cuntz-Li-algebra}.) We show that $\mathfrak{A}[N \rtimes H,M]$ is generated by an inverse semigroup of partial isometries denoted by $T$. We show that $\mathfrak{A}[N \rtimes H,M]$ is isomorphic to the $C^{*}$-algebra of the groupoid $\mathcal{G}_{tight}$, considered in \cite{Ex}, of the inverse semigroup $T$. We also identify the groupoid $\mathcal{G}_{tight}$ explicitly and show that $\mathcal{G}_{tight}$ is isomorphic to $\overline{N}\rtimes G|_{\overline{M}}$.  The author had done a similar analysis for the Cuntz-Li algebra associated to the ring $\mathbb{Z}$ in \cite{Sundar1}. At the end of this paper, we prove a duality result  analogous to the duality result obtained in \cite{Cuntz-Li-1}.
\section{Semidirect products and the Cuntz-Li relations}
 Let $G=N \rtimes H$ be a semidirect product and let $M$ be a normal subgroup of $N$. Let $P:=\{a \in H:aMa^{-1} \subset M\}$. Then $P$ is a semigroup containing the identity $e$.  Assume that the following holds.
\begin{enumerate}
 \item[(C1)] The group $H=PP^{-1}=P^{-1}P$.
  \item[(C2)] For every $a \in P$, the subgroup $aMa^{-1}$ is of finite index in $M$.
  \item[(C3)] The intersection $\displaystyle \bigcap_{a \in P}aMa^{-1} = \{e\}$ where $e$ denotes the identity element of $G$.
\end{enumerate}

Let $\mathcal{U}=\{aMa^{-1}:a \in H \}$. In \cite{Quigg-Landstad}, the following conditions were required to be satisfied. (Cf. Section 2 in \cite{Quigg-Landstad}.)
\begin{enumerate}
 \item[(E1)] Given $U,V \in \mathcal{U}$, there exists $W \in \mathcal{U}$ such that $W \subset U \cap V$.
 \item[(E2)] If $U,V \in \mathcal{U}$ and $U \subset V$ then $U$ is of finite index in $V$.
 \item[(E3)] The intersection $\displaystyle \bigcap_{U \in \mathcal{U}}U=\{e\}$.
\end{enumerate}

 We claim that (E1) is equivalent to the condition $H=PP^{-1}$. Assume (E1). Let $a \in H$ be given. Then there exists $c \in H$ such that $a^{-1}Ma \cap M \supseteq cMc^{-1}$. Then $c \in P$ and $ac \in P$. Note that $a=(ac)c^{-1} \in PP^{-1}$. Thus we have $H=PP^{-1}$.

Now suppose $H=PP^{-1}$. First note that for every $a,b \in P$, $aP \cap bP$ is non-empty. Now let $c,d \in H$ be given. Write $c=a_{1}a_{2}^{-1}$ and $d=b_{1}b_{2}^{-1}$ with $a_{i},b_{i} \in P$. Choose $\alpha, \beta \in P$ such that $a_{1}\alpha=b_{1}\beta$. Let $a:=a_{1}\alpha$. Then $c^{-1}a=a_{2}\alpha \in P$. Similarly $d^{-1}a \in P$. Hence $aMa^{-1} \subset cMc^{-1} \cap dMd^{-1}$. Thus (E1) holds.

Given (E1), note that (E3) is equivalent to (C3). For if $a \in H$, there exists $b \in P$ such that $aMa^{-1} \cap M \supseteq bMb^{-1}$. Thus for every $a \in H$, $aMa^{-1} \supseteq \bigcap_{b \in P}bMb^{-1}$. Hence $\bigcap_{U \in \mathcal{U}}U=\bigcap_{a \in P}aMa^{-1}$. Thus given (E1), (E3) is equivalent to (C3). Clearly (E2) is equivalent to (C2).

\begin{rmrk}
 In \cite{Quigg-Landstad},  the Cuntz-Li algebra associated to the pair ( Cf. defn \ref{Cuntz-Li-algebra}) $(N \rtimes H,M)$ was considered when $H$ is abelian. (Cf. Hypothesis 9.2 and Theorem 9.11 in \cite{Quigg-Landstad}.) Here, we consider a slightly more general situation. We assume $H=P^{-1}P=PP^{-1}$.
\end{rmrk}

\begin{rmrk}
 The condition $H=P^{-1}P=PP^{-1}$ is equivalent to saying that $P$ generates $H$ and $P$ is right and left reversible i.e. given $a,b \in P$, the intersections $Pa\cap Pb$ and $aP\cap bP$ are non-empty. Cancellative semigroups which are right (or left) reversible are called Ore semigroups. For more details on Ore semigroups, we refer to \cite{Clifford}.

 A semigroup $P$ is called right reversible (left reversible) if $Pa \cap Pb$ (if $aP \cap bP$) is non-empty for every $a,b \in P$.
\end{rmrk}

Throughout this article, whenever we write $G=N\rtimes H$ and $M$ is a normal subgroup of $N$, we assume that conditions (C1), (C2) and (C3) hold. For $a \in P$, let $M_{a}=aMa^{-1}$. We will use this notation throughout.

\begin{lmma}
\label{The subgroup N}
 Let $G=N\rtimes H$ and $M$ be a normal subgroup of $N$. Let $N_{0}:=\displaystyle \bigcup_{a \in P}a^{-1}Ma$. Then $N_{0}$ is a subgroup of $N$ and is invariant under conjugation by $H$.
\end{lmma}
\textit{Proof.} First observe that $N_{0}$ is closed under inversion. Let $a,b \in P$ be given. Choose an element $c$ in the intersection $Pa \cap Pb$. Then $a^{-1}Ma \subset c^{-1}Mc$ and $b^{-1}Mb \subset c^{-1}Mc$. Now it follows that $N_{0}$ is closed under multiplication. Thus $N_{0}$ is a subgroup of $N$.

Obviously $N_{0}$ is invariant under conjugation by $P^{-1}$. Let $a,b \in P$ be given. Since $P$ is right reversible, there exists $c,d \in P$ such that $ab^{-1}=c^{-1}d$. Now observe that $a(b^{-1}Mb)a^{-1}=c^{-1}(dMd^{-1})c \subset c^{-1}Mc$
  Thus it follows that $N_{0}$ is closed under conjugation by $P$. This completes the proof. \hfill $\Box$
\begin{rmrk}
\label{key remark}
 As a consequence of Lemma \ref{The subgroup N}, we may very well assume as in \cite{Quigg-Landstad} that $N=\displaystyle \bigcup_{a \in P}a^{-1}Ma$. 
\end{rmrk}
Let us consider a few examples which fits the setup that we are considering.
\begin{xmpl}[\cite{Cuntz-Li}]
\label{Cuntz-Li}
 Let $R$ be an integral domain such that for every non-zero $m \in R$, the ideal generated by $m$ is of finite index in $R$. Assume that $R$ is not a field. We denote the field of fractions of $R$ by $Q$ and the set of non-zero elements in $Q$   by $Q^{\times}$. The multiplicative group $Q^{\times}$ acts on $Q$ by multiplication. Now let $N:=Q$, $H:=Q^{\times}$ and $M:=R$. Then $P=R^{\times}$ where $R^{\times}$ denotes the set of non-zero elements in $R$. Then conditions (C1)-(C3) hold for the pair $(N \rtimes H, M)$. 
\end{xmpl}
\begin{xmpl}[\cite{Quigg-Landstad}]
 Let $F$ be a finite group and consider the direct sum $N:=\oplus_{\mathbb{Z}}F$. Then $H:=\mathbb{Z}$ acts on $N$ by shifting. Let $M:=\oplus_{\mathbb{N}}F$ be the normal subgroup of $N$. Then it is easily verifiable that the pair $(N \rtimes H,M)$ satisfies the hypothesis (C1)-(C3).
\end{xmpl}
In the following two examples, we think of elements of $\mathbb{Q}^{n}$ as column vectors.
\begin{xmpl}
\label{Exel}
 Let $A$ be a $n \times n$ integer dilation matrix. In other words, $A$ is an $n\times n$ matrix with integer entries such that every complex eigen value of $A$ has absolute value greater than $1$. Note that $A$ is invertible over $\mathbb{Q}$ and $|\det(A)|>1$. The matrix $A$ acts on $\mathbb{Q}^{n}$ by matrix multiplication and thus induces an action of $\mathbb{Z}$ on $\mathbb{Q}^{n}$. We let the generator $1$ of $\mathbb{Z}$ act on $\mathbb{Q}^{n}$ by $1.v=Av$ for $v \in \mathbb{Q}^{n}$. Let $N:=\mathbb{Q}^{n}$, $H:=\mathbb{Z}$ and $M:=\mathbb{Z}^{n}$. Then $P=\mathbb{N}$. Let us verify the hypothesis (C1)-(C3).
\begin{enumerate}
 \item[(C1)] Note that $H$ is abelian and $H=PP^{-1}=P^{-1}P$.
 \item[(C2)] For $r \geq 0$, the index of $A^{r}\mathbb{Z}^{n}$ is of finite index in $\mathbb{Z}^{n}$ and in fact its index is $|\det(A)|^{r}$.
 \item[(C3)] Lemma 4.1 of \cite{Exel-Rae} implies that the operator norm $||A^{-m}||$ converges to $0$ as $m$ tends to infinity. Thus if $\displaystyle 0 \neq v \in \bigcap_{r=0}^{\infty}A^{r}\mathbb{Z}^{n}$, then for every $m \geq 0$, $A^{-m}v \in \mathbb{Z}^{n}$. Thus we have $1 \leq ||A^{-m}v|| \leq ||A^{-m}||||v||$ which is a contradiction. Thus (C3) holds.
\end{enumerate}
\end{xmpl}
The case $n=1$ and $A=p$ where $p$ is a prime number was discussed in \cite{Larsen-Li-1}.
In the previous example, we can consider integer matrices other than dilation matrices. It is possible that (C3) is satisfied for an integer matrix $A$ such that $|\det(A)|>1$ and $\bigcap_{r>0}A^{r}\mathbb{Z}^{n}=\{0\}$ without $A$ being a dilation matrix. In fact we have the following nice characterisation of condition (C3) when $n=2$.
\begin{lmma}
 Let $A$ be a $2\times 2$ matrix with integer entries. Assume that $|\det(A)|>1$.  Then the following are equivalent.
\begin{enumerate}
 \item The intersection $\displaystyle \bigcap_{r \geq 0}A^{r}\mathbb{Z}^{2}$ is trivial.
 \item Neither $1$ nor $-1$ is an eigen value of $A$.
\end{enumerate}
\end{lmma}
\textit{Proof.} Suppose $\bigcap_{r \geq 0}A^{r}\mathbb{Z}^{2}=\{0\}$. If $1$ is an eigen value of $A$ then there exists a non-zero $v \in \mathbb{Q}^{2}$ such that $Av=v$. By clearing denominators, we can assume that $v \in \mathbb{Z}^{2}$. Then clearly $v \in \displaystyle \bigcap_{r \geq 0}A^{r}\mathbb{Z}^{2}$. Thus we have shown that $1$ is not an eigen value of $A$. Similarly we can show $-1$ is not an eigen value of $A$.

Now assume that neither $1$ nor $-1$ is an eigen value of $A$. Let $\Gamma_{r}:=A^{r}\mathbb{Z}^{2}$ and $\Gamma:= \bigcap_{r\geq 0}\Gamma_{r}$. Since $\Gamma \subset \Gamma_{r} \subset \mathbb{Z}^{2}$, we have  $[\mathbb{Z}^{2}:\Gamma] \geq [\mathbb{Z}^{2}:\Gamma_{r}]=|\det(A)|^{r}$. Hence $\Gamma$ cannot be of finite index in $\mathbb{Z}^{2}$. This implies that $\Gamma$ is of rank atmost $1$. If $\Gamma$ is rank $1$ then there exists a non-zero $v \in \mathbb{Z}^{2}$ such that $\Gamma=\mathbb{Z}v$. But $A:\Gamma \to \Gamma$ is a bijection. Thus it must  either be multiplication by $1$ or by $-1$. In other words, $v$ is an eigen vector for $A$ with eigen value $1$ or $-1$. This is a contradiction. Thus $\Gamma$ cannot be of rank $1$ which in turn implies $\Gamma=\{0\}$. This completes the proof. \hfill $\Box$.

The matrix $A:= \begin{bmatrix}
             0 & 2 \\
             1 & -2
            \end{bmatrix}$ has eigen values $\sqrt{3}-1$ and $-\sqrt{3}-1$. But $A$ is not a dilation matrix but still (C3) holds for $A$.
\begin{rmrk}
 It is not clear to the author whether (C3) can be characterised in terms of eigen values of the matrix in the higher dimensional case.
\end{rmrk}
Let us now consider an example where $H$ is non-abelian.
\begin{xmpl}
\label{matrix group}
 Let $N=\mathbb{Q}^{n}$ and $H$ be a subgroup of $GL_{n}(\mathbb{Q})$ containing the non-zero scalars. Just as in Example \ref{Exel}, $H$ acts on $N$ by matrix multiplication. Let $M=\mathbb{Z}^{n}$. Then $P$ consists of elements of $H$ whose entries are integers. 
\begin{enumerate}
 \item[(C1)] Let $A \in H$ be given. Then there exists a non-zero integer $m$ such that $mA=Am \in P$. Hence $H=PP^{-1}=P^{-1}P$.
 \item[(C2)] For $A \in P$, the subgroup $A\mathbb{Z}^{n}$ is of finite index and its index is $|\det(A)|$.
 \item[(C3)] Since $\displaystyle \bigcap_{ m \in \mathbb{Z}^{\times}}m\mathbb{Z}^{n}=\{0\}$, it follows that $\displaystyle \bigcap_{A \in P}A\mathbb{Z}^{n}=\{0\}$.
\end{enumerate}
\end{xmpl}
\begin{dfn}
\label{Cuntz-Li-algebra}
Let $G:=N\rtimes H$ be a semidirect product and $M$ be a normal subgroup of $N$ such that (C1)-(C3) holds. We let $\mathfrak{A}[N \rtimes H,M]$ be the universal $C^{*}$-algebra generated by a set of isometries $\{s_{a}:a \in P\}$ and a set of unitaries $\{u(m): m \in M\}$ satisfying the following relations.
\begin{equation*}
 \begin{split}
  s_{a}s_{b}&=s_{ab}\\
  u(m)u(n)&=u(mn)\\
  s_{a}u(m)&=u(ama^{-1})s_{a}\\
  \displaystyle \sum_{k \in M/M_{a}}u(k)e_{a}u(k)^{-1}&=1
 \end{split}
\end{equation*}
where $e_{a}$ denotes the final projection of $s_{a}$.
\end{dfn}
Note that $u(k)e_{a}u(k)^{-1}$ depends only on the coset $k(M_{a})$. Moreover if $k_{1}$ and $k_{2}$ lie in different cosets of $M_{a}$ then $u(k_{1})e_{a}u(k_{1})^{-1}$ and $u(k_{2})e_{a}u(k_{2})^{-1}$ are orthogonal.

For $a\in P$ and $m \in M$, consider the operators $S_{a}$ and $U(m)$  on $\ell^{2}(M)\otimes \ell^{2}(H)$ defined as follows
\begin{equation*}
 \begin{split}
  S_{a}(\delta_{n}\otimes \delta_{b}):&=\delta_{ana^{-1}}\otimes \delta_{ab}\\
  U(m)(\delta_{n}\otimes \delta_{b}):&=\delta_{mn}\otimes \delta_{b}.
 \end{split}
\end{equation*}
Then $s_{a} \to S_{a}$ and $u(m)\to U(m)$ gives a representation of $\mathfrak{A}[N \rtimes H,M]$ on the Hilbert space $\ell^{2}(M)\otimes \ell^{2}(H)$. Let us call this representation the regular representation and denote its image by $\mathfrak{A}_{r}[N \rtimes H,M]$.

\begin{rmrk}
 It should be noted that the regular representation for integral domains considered in \cite{Cuntz-Li} is different from ours.
\end{rmrk}

\section{An Inverse semigroup for the Cuntz-Li relations}
The main aim of this section is to show that the $C^{*}$-algebra $\mathfrak{A}[N \rtimes H,M]$ is generated by an inverse semigroup of partial isometries.  We begin with a lemma similar to  Lemma 1 of Section 3.1 in \cite{Cuntz-Li}.

\begin{lmma}
\label{decomposition}
 For every $a,b \in P$, one has \[
e_{a}= \displaystyle \sum_{k \in M/M_{b}}u(aka^{-1})e_{ab}u(aka^{-1})^{-1}                                     
                                   \]
\end{lmma}
\textit{Proof.} One has 
\begin{equation*}
 \begin{split}
  e_{a}&= s_{a}s_{a}^{*} \\
       &=s_{a}\big(\sum_{k \in M/M_{b}}u(k)e_{b}u(k)^{-1}\big)s_{a}^{*} \\
       &=\sum_{k \in M/M_{b}}u(aka^{-1})s_{a}e_{b}s_{a}^{*}u(aka^{-1})^{-1} \\
       &= \sum_{k \in M/M_{b}}u(aka^{-1})e_{ab}u(aka^{-1})^{-1}
 \end{split}
\end{equation*}
This completes the proof. \hfill $\Box$

Let $X$ be the linear span of $\{u(k)e_{b}u(k)^{-1}:b \in P, k \in M\}$. Denote the set of projections  in $X$ by $F$. By Lemma \ref{decomposition} and the left reversibility of $P$, it follows that $f \in F$ if and only if there exists $b \in P$ such that $f$ is in the linear span of $\{u(k)e_{b}u(k)^{-1}\}$. The following lemma is an immediate corollary of Lemma \ref{decomposition} and the fact that $P$ is left reversible.

\begin{lmma}
\label{semigroup of projections}
 The set $F$ is a commutative semigroup of projections. Moreover $F$ is invariant under the maps $x \to s_{b}xs_{b}^{*}$ for every $b \in P$ and $x \to u(m)xu(m)^{-1}$ for every $m \in M$.
\end{lmma}

Now we show that $F$ is also invariant under conjugation by $s_{a}^{*}$ for every $a \in P$.

\begin{lmma}
\label{conjugation}
 Let $a \in P$ be given. If $f \in F$, then $s_{a}^{*}fs_{a} \in F$. Moreover, $s_{a}^{*}u(m)e_{b}u(m)^{-1}s_{a}$ is in the linear span of $\{u(k)e_{a^{-1}c}u(k)^{-1} \}$ where $c$ is any element in $aP \cap bP$.
\end{lmma}
\textit{Proof.} Let $a \in P$ and $f \in F$ be given. First observe that $s_{a}^{*}fs_{a}$ is selfadjoint. Also 
\begin{equation*}
\begin{split}
 (s_{a}^{*}fs_{a})^{2}&=s_{a}^{*}fs_{a}s_{a}^{*}fs_{a} \\
                      &=s_{a}^{*}fe_{a}fs_{a} \\
                      &=s_{a}^{*}e_{a}fs_{a} \mbox{ ( Since $F$ is commutative ) }\\
                      &=s_{a}^{*}fs_{a}
\end{split}
\end{equation*}
Thus $s_{a}^{*}fs_{a}$ is a projection. Now to show that $s_{a}^{*}fs_{a} \in F$, it is enough to consider the case when $f=u(m)e_{b}u(m)^{-1}$. Now let $c \in aP \cap bP$ and write $c=a\alpha=b \beta$ with $\alpha,\beta \in P$.

Let $r_{1},r_{2},\cdots, r_{n}$ be distinct representatives of $M/M_{\beta}$. Then by Lemma \ref{decomposition}, it follows that 
\begin{align*}
 s_{a}^{*}u(m)e_{b}u(m)^{-1}s_{a}&= \displaystyle \sum_{i=1}^{n}s_{a}^{*}u(mbr_{i}b^{-1})e_{b\beta}u(mbr_{i}b^{-1})^{-1}s_{a} \\
                                 &= \displaystyle \sum_{i=1}^{n}s_{a}^{*}u(mbr_{i}b^{-1})e_{a\alpha}u(mbr_{i}b^{-1})^{-1}s_{a}
\end{align*}
The term $s_{a}^{*}u(mbr_{i}b^{-1})e_{a\alpha}u(mbr_{i}b^{-1})^{-1}s_{a}$ survives if and only if $e_{a\alpha}u(mbr_{i}b^{-1})s_{a} \neq 0$ and that is if and only if $e_{a\alpha}u(mbr_{i}b^{-1})e_{a}u(mbr_{i}b^{-1})^{-1} \neq 0$. But by Lemma \ref{decomposition} this happens precisely when there exists $t_{i} \in M/M_{\alpha}$ such that $mbr_{i}b^{-1} \equiv at_{i}a^{-1} mod M_{a\alpha}$.

Let \[
A:=\{ i : \mbox{ There exists $t_{i}$ such that $mbr_{i}b^{-1} \equiv at_{i}a^{-1} \mod M_{a\alpha}$} \}.     
    \]
 For every $i \in A$, choose $t_{i}$ such that $mbr_{i}b^{-1} \equiv at_{i}a^{-1} \mod M_{a\alpha}$. Now we have
\begin{equation*}
\begin{split}
 s_{a}^{*}u(m)e_{b}u(m)^{-1}s_{a}&= \displaystyle \sum_{i=1}^{n}s_{a}^{*}u(mbr_{i}b^{-1})e_{a\alpha}u(mbr_{i}b^{-1})^{-1}s_{a}\\
                                 &=\displaystyle \sum_{i \in A}s_{a}^{*}u(mbr_{i}b^{-1})e_{a\alpha}u(mbr_{i}b^{-1})^{-1}s_{a}\\
                                 &=\displaystyle \sum_{i \in A} s_{a}^{*}u(at_{i}a^{-1})e_{a\alpha}u(at_{i}a^{-1})^{-1}s_{a}\\
                                 &=\displaystyle \sum_{i \in A} u(t_{i})s_{a}^{*}e_{a\alpha}s_{a}u(t_{i})^{-1}\\
                                 &=\displaystyle \sum_{i \in A}u(t_{i})e_{\alpha}u(t_{i})^{-1} 
\end{split}
\end{equation*}
This completes the proof. \hfill $\Box$

Let us isolate the computation in the previous lemma in a remark. This will be used later.
\begin{rmrk}
\label{conjugation by s}
 Let $a,b \in P$ be given. Let $c \in aP \cap bP$. Choose $\alpha$ and $\beta$ in $P$ such that $c=a\alpha=b\beta$. Conjugation by $a$ sends $M_{\alpha}$ to $M_{c}$. Thus we get a map denoted $\pi_{\alpha}^{a}:M/M_{\alpha} \to M/M_{c}$. Similarly conjugation by $b$ gives a map $\pi_{\beta}^{b}:M/M_{\beta} \to M/M_{c}$. Note that both $\pi_{\alpha}^{a}$ and $\pi_{\beta}^{b}$ are injective. Denote the quotient map $M \to M/M_{c}$ by $q_{c}$. For $m \in M$, define
\[
 A_{m}:=\{r \in M/M_{\beta}: q_{c}(m)\pi_{\beta}^{b}(r) \in \pi_{\alpha}^{a}(M/M_{\alpha}) \}.
\]
Then the computation in Lemma \ref{conjugation} can be restated as follows
\begin{equation*}
 s_{a}^{*}u(m)e_{b}u(m)^{-1}s_{a}=\displaystyle \sum_{r \in A_{m}}u\Big(\big(\pi^{a}_{\alpha}\big)^{-1}(q_{c}(m)\pi^{b}_{\beta}(r))\Big)e_{\alpha}u\Big(\big(\pi^{a}_{\alpha}\big)^{-1}(q_{c}(m)\pi^{b}_{\beta}(r))\Big)^{-1}.
\end{equation*}
\end{rmrk}
Now we show that $\mathfrak{A}[N \rtimes H,M]$ is generated by an inverse semigroup of partial isometries.

\begin{ppsn}
\label{Inverse semigroup}
 Let $T:=\{s_{a}^{*}u(m)fu(m^{'})s_{a^{'}}: m,m^{'} \in M, a,a^{'} \in P, \mbox{and~}f \in F \}$. Then $T$ is an inverse semigroup of partial isometries containing $0$. Moreover the set of projections in $T$ coincides exactly with $F$. Also the linear span of $T$ is a dense $*$-subalgebra of $\mathfrak{A}[N \rtimes H,M]$.
\end{ppsn}
\textit{Proof.} The fact that $T$ is closed under multiplication follows from the following calculation. Let $a_{1},a_{2},b_{1},b_{2} \in P$, $m_{1},m_{2},n_{1},n_{2} \in M$ and $e,f \in F$ be given. Choose $c \in Pb_{1} \cap Pa_{2}$ and write $c$ as $c=\beta b_{1}= \alpha a_{2}$. Observe that
\begin{equation*}
 \begin{split}
&s_{a_{1}}^{*}u(m_{1})eu(m_{2})s_{a_{2}}s_{b_{1}}^{*}u(n_{1})fu(n_{2})s_{b_{2}}\\
 & = s_{a_{1}}^{*}u(m_{1}m_{2})u(m_{2}^{-1})eu(m_{2})s_{\alpha}^{*}s_{\alpha}s_{a_{2}}s_{b_{1}}^{*}s_{\beta}^{*}s_{\beta}u(n_{1})fu(n_{1}^{-1})u(n_{1}n_{2})s_{b_{2}} \\
 &=s_{a_{1}}^{*}u(m_{1}m_{2})u(m_{2}^{-1})eu(m_{2})s_{\alpha}^{*}s_{\alpha a_{2}}s_{\beta b_{1}}^{*}s_{\beta}u(n_{1})fu(n_{1}^{-1})u(n_{1}n_{2})s_{b_{2}} \\
 &=s_{a_{1}}^{*}u(m_{1}m_{2})s_{\alpha}^{*}s_{\alpha}u(m_{2}^{-1})eu(m_{2})s_{\alpha}^{*}s_{c}s_{c}^{*}s_{\beta}u(n_{1})fu(n_{1}^{-1})s_{\beta}^{*}s_{\beta}u(n_{1}n_{2})s_{b_{2}}\\
 &=s_{a_{1}}^{*}s_{\alpha}^{*}u(\alpha m_{1}m_{2} \alpha^{-1})\big(s_{\alpha}u(m_{2}^{-1})eu(m_{2})s_{\alpha}^{*}\big)e_{c}\big(s_{\beta}u(n_{1})fu(n_{1}^{-1})s_{\beta}^{*}\big)u(\beta n_{1}n_{2} \beta^{-1})s_{\beta}s_{b_{2}} \\
 &=s_{\alpha a_{1}}^{*}u(\alpha m_{1}m_{2} \alpha^{-1})\big(s_{\alpha}u(m_{2}^{-1})eu(m_{2})s_{\alpha}^{*}\big)e_{c}\big(s_{\beta}u(n_{1})fu(n_{1}^{-1})s_{\beta}^{*}\big)u(\beta n_{1}n_{2} \beta^{-1})s_{\beta b_{2}}\\
&=s_{\alpha a_{1}}^{*}u(\alpha m_{1}m_{2} \alpha^{-1})(s_{\alpha}\tilde{e}s_{\alpha}^{*})e_{c}(s_{\beta}\tilde{f}s_{\beta}^{*})u(\beta n_{1}n_{2} \beta^{-1})s_{\beta b_{2}} \\
\end{split}
\end{equation*}
where $\tilde{e}=u(m_{2}^{-1})eu(m_{2})$ and $\tilde{f}=u(n_{1})fu(n_{1})^{-1}$. The above calculation together with Lemma \ref{semigroup of projections} implies that $T$ is closed under multiplication. Obviously $T$ is closed under the involution $*$.

Now let us show that every element of $T$ is a partial isometry. Let $v:=s_{a}^{*}u(m)fu(m^{'})s_{a^{'}}$ be an element of $T$. Then 
\[
vv^{*}=s_{a}^{*}\Big (u(m)\big (fu(m^{'})e_{a^{'}}u(m^{'})^{-1}f\big)u(m)^{-1}\Big)s_{a}\]

Now Lemma \ref{semigroup of projections} and Lemma \ref{conjugation} implies that $vv^{*} \in F$. Thus we have shown that every element of $T$ is a partial isometry and the set of projections in $T$ coincides with $F$. In other words $T$ is an inverse semigroup.

Since $T$ is  closed under multiplication and involution, it follows that the linear span of $T$ is a $*$-algebra.
Moreover $T$ contains $\{s_{a}: a \in P\}$ and $\{u(m):m \in M\}$. Thus the linear span of $T$ is dense in $\mathfrak{A}[N \rtimes H,M]$. This completes the proof. \hfill $\Box$

The following equality will be used later. Let $a_{1},a_{2},b_{1},b_{2} \in P$  and $m_{1},m_{2} \in M$ be given. Choose $c \in Pb_{1} \cap Pa_{2}$ and write $c$ as $c=\beta b_{1} = \alpha a_{2}$.  Now the computation in Proposition \ref{Inverse semigroup} gives the following equality
\begin{equation}
\label{commutation rule}
 s_{a_{1}}^{*}u(m_{1})s_{b_{1}}s_{a_{2}}^{*}u(m_{2})s_{b_{2}} = s_{\beta a_{1}}^{*}u(\beta m_{1} \beta^{-1})e_{c}u(\alpha m_{2} \alpha^{-1})s_{\alpha b_{2}}
\end{equation}

\begin{rmrk}
\label{good remark}
 We also need the following fact. If $v \in T$, let us denote its image in the regular representation by $V$. Observe that $v \neq 0$ if and only if $V \neq 0$. This is clear for projections in $T$. Now let $v \in T$ be a non-zero element. Then $vv^{*} \in F$ is non-zero. Thus $VV^{*} \neq 0$ which implies $V\neq 0$.
\end{rmrk}

In the remainder of this article, we reserve the letter $T$ to denote the inverse semigroup in Proposition \ref{Inverse semigroup} and $F$ to denote the set of projections in $T$.

\section{Tight representations of  inverse semigroups}

In this section, we  show that the identity representation of $T$ in $\mathfrak{A}[N \rtimes H,M]$ is tight in the sense of Exel and the $C^{*}$-algebra of the tight groupoid associated to $T$ is isomorphic to $\mathfrak{A}[N \rtimes H,M]$. First let us recall the notion of tight characters and tight representations from \cite{Ex}.
\begin{dfn}
 Let $S$ be an inverse semigroup with $0$. Denote the set of projections in $S$ by $E$. A character for $E$ is a map $x:E \to \{0,1\}$ such that 
\begin{enumerate}
 \item the map $x$ is a semigroup homomorphism, and
 \item $x(0)=0$.
\end{enumerate}
\end{dfn}
We denote the set of characters of $E$ by $\widehat{E_{0}}$. We consider $\widehat{E_{0}}$ as a locally compact Hausdorff topological space where the topology on $\widehat{E_{0}}$ is the subspace topology induced from the product topology on 
$\{0,1\}^{E}$. 

For a character $x$ of $E$, let $A_{x}:=\{e \in E: x(e)=1\}$. Then $A_{x}$ is a nonempty set satisfying the following properties.
\begin{enumerate}
 \item[(1)]The element $0 \notin A_{x}$.
  \item[(2)]If $e \in A_{x}$ and $f \geq e$ then $f \in A_{x}$.
  \item[(3)]If $e,f \in A_{x}$ then $ef \in A_{x}$.
 \end{enumerate}

Any nonempty subset $A$ of $E$ for which $(1),(2)$ and $(3)$ are satisfied is called a filter. Moreover if $A$ is a filter then the indicator function $1_{A}$ is a character. Thus there is a bijective correspondence between the set of characters and filters. A filter is called an ultrafilter if it is maximal. We also call a character $x$ maximal or an ultrafilter if its support $A_{x}$ is maximal. The set of maximal characters is denoted by $\widehat{E_{\infty}}$ and its closure in $\widehat{E_{0}}$ is denoted by $\widehat{E_{tight}}$.

We refer to \cite{Sundar1} (Corollary 3.3) for the proof of the following lemma.
\begin{lmma}
\label{maximality}
 Let $A$ be a unital $C^{*}$-algebra and $E \subset A$ be an inverse semigroup of projections containing $\{0,1\}$. Suppose that $E$ contains a finite set $\{e_{1},e_{2},\cdots,e_{n}\}$ of mutually orthogonal projections such that $\sum_{i=1}^{n}e_{i}=1$. Then for every maximal character $x$ of $E$, there exists a unique $e_{i}$ for which $x(e_{i})=1$.
\end{lmma}

Let us recall the notion of tight representations of semilattices from \cite{Ex} and from \cite{Ex1}. The only semilattice we consider is that of an inverse semigroup of projections or in other words the idempotent semilattice of an inverse semigroup. Also our semilattice contains a maximal element $1$. First let us recall the notion of a cover from \cite{Ex}.
\begin{dfn}
 Let $E$ be an inverse semigroup of projections containing $\{0,1\}$ and $Z$ be a subset of $E$. A subset $F$ of $Z$ is called a cover for $Z$ if given a non-zero element $z \in Z$ there exists an $f \in F$ such that $fz\neq 0$. A cover $F$ of $Z$ is called a finite cover if $F$ is finite.
\end{dfn}

The following definition is actually Proposition 11.8 in \cite{Ex}
\begin{dfn}
 Let $E$ be an inverse semigroup of projections containing $\{0,1\}$. A representation $\sigma:E \to \mathcal{B}$ of the semilattice $E$ in a Boolean algebra $\mathcal{B}$ is said to be tight if $\sigma(0)=0$ and given $e\neq 0$ in $E$ and for every finite cover $F$ of the interval $[0,e]:=\{x \in E: x \leq e\}$, one has $\sup_{f \in F} \sigma(f)=\sigma(e)$.
\end{dfn}

Let $A$ be a unital $C^{*}$ algebra and $S$ be an inverse semigroup containing $\{0,1\}$. Denote the set of projections in $S$ by $E$. Let $\sigma:S \to A$ be a unital representation of $S$ as partial isometries in $A$. Let $\sigma(C^{*}(E))$ be the $C^{*}-$subalgebra in $A$ generated by $\sigma(E)$. Then $\sigma(C^{*}(E))$ is a unital, commutative $C^{*}-$algebra and hence the set of projections in it is a Boolean algebra which we denote by $\mathcal{B}_{\sigma(C^{*}(E))}$. We say the representation $\sigma$ is \textbf{tight} if the representation $\sigma:E \to \mathcal{B}_{\sigma(C^{*}(E))}$ is \textbf{tight}. The proof of the following lemma can be found in \cite{Sundar1} (Lemma 3.6, page 7).
\begin{lmma}
\label{tightness}
 Let $X$ be a compact metric space and  $E \subset C(X)$ be an inverse semigroup
of projections containing $\{0,1\}$. Suppose that for every finite set of
projections $\{f_{1},f_{2},\cdots,f_{m}\}$ in $E$, there exists a finite set of
mutually orthogonal non-zero projections $\{e_{1},e_{2},\cdots,e_{n}\}$ in $E$ and a 
matrix $(a_{ij})$ such that
\begin{align*}
\sum_{i=1}^{n}e_{i}&=1 \\
f_{i}&=\sum_{j}a_{ij}e_{j}. 
\end{align*}
Then the identity representation of $E$ in $C(X)$ is tight.
\end{lmma}

As in \cite{Sundar1}, we prove that the identity representation of $T$ in $\mathfrak{A}[N \rtimes H,M]$ is tight.
\begin{ppsn}
 \label{tightness of the identity representation}
The identity representation of $T$ in $\mathfrak{A}[N \rtimes H,M]$ is tight.
\end{ppsn}
\textit{Proof.} We apply Lemma \ref{tightness}. Let $\{f_{1},f_{2},\cdots,f_{n}\}$ be a finite set of projections in $T$. By definition, given $i$ there exists $a_{i} \in P$ such that $f_{i}$ is in the linear span of $\{u(k)e_{a_{i}}u(k)^{-1}\}$. Let $c \in \displaystyle \bigcap_{i=1}^{n} a_{i}P$. By Lemma \ref{decomposition}, it follows that for every $i$, $f_{i}$ is in the linear span of $\{u(k)e_{c}u(k)^{-1}: k \in M/cMc^{-1}\}$. Appealing to Lemma \ref{tightness}, we can conclude that the identity representation of $T$ in $\mathfrak{A}[N \rtimes H,M]$ is tight. This completes the proof. \hfill $\Box$

Now we show that $\mathfrak{A}[N \rtimes H,M]$ is isomorphic to the $C^{*}$-algebra of the groupoid $\mathcal{G}_{tight}$ associated to $T$. For the convenience of the reader, we recall the construction of the groupoid $\mathcal{G}_{tight}$, considered in \cite{Ex}, associated to an inverse semigroup with $0$. 

Let $S$ be an inverse semigroup with $0$ and let $E$ denote its set of projections. Note that $S$ acts on $\widehat{E_{0}}$ partially. 
For $x \in \widehat{E_{0}}$ and  $s \in S$, define $(x.s)(e)=x(ses^*)$. Then
\begin{itemize}
 \item The map $x.s$ is a semigroup homomorphism, and 
 \item $(x.s)(0)=0$.
\end{itemize}
But $x.s$ is nonzero if and only if $x(ss^*)=1$.  For $s \in S$, define the domain
and range of $s$ as 
\begin{align*}
 D_s:&=\{x \in \widehat{E_{0}}: x(ss^*)=1\} \\
 R_s:&=\{x \in \widehat{E_{0}}: x(s^*s)=1\} 
\end{align*}
Note that both $D_s$ and $R_s$ are compact and open. Moreover $s$ defines a
homeomorphism from $D_s$ to $R_s$ with $s^*$ as its inverse. Also observe that
$\widehat{E_{tight}}$ is invariant under the action of $S$.

Consider the transformation groupoid $\Sigma:=\{(x,s):x\in D_s\}$ with the
composition and the inversion being given by:
 \begin{align*}
  (x,s)(y,t):&=(x,st) \textrm{~if~} y=x.s\\
   (x,s)^{-1}:&=(x.s,s^*) 
 \end{align*}
Define an equivalence relation $\sim$ on $\Sigma$ as $(x,s)\sim (y,t)$ if $x=y$
and if there exists an $e \in E$ such that $x \in D_e$ for which $es=et$. Let
$\mathcal{G}=\Sigma/\sim$. Then $\mathcal{G}$ is a groupoid as the product and
the inversion respects the equivalence relation $\sim$. Now we describe a
topology on $\mathcal{G}$ which makes $\mathcal{G}$ into a topological groupoid.

For $s \in S$ and $U$ an open subset of $D_s$, let $\theta(s,U):=\{[x,s]: x \in
U \}$. We refer to \cite{Ex} for the proof of the following  proposition. We
denote $\theta(s,D_s)$ by $\theta_s$. 
\begin{ppsn}
 The collection $\{\theta(s,U): s \in S, U \textrm{~open in ~} D_s \}$ forms a
basis for a topology on $\mathcal{G}$. The groupoid $\mathcal{G}$ with this
topology is a topological groupoid whose unit space can be identified with
$\widehat{E_{0}}$.
Also one has the following.
\begin{enumerate}
 \item For $s,t \in S$, $\theta_s \theta_t=\theta_{st}$,
 \item For $s \in S$, $\theta_{s}^{-1}=\theta_{s^{*}}$, 
 \item For $s \in S$, $\theta_{s}$  is compact, open and Hausdorff, and
 \item The set $\{1_{\theta_s}: s \in T\}$ generates the $C^*$-algebra
$C^*(\mathcal{G})$.
\end{enumerate}
\end{ppsn}

We define the groupoid $\mathcal{G}_{tight}$ to be the reduction of the groupoid
$\mathcal{G}$ to $\widehat{E_{tight}}$. In \cite{Ex}, it is shown that the representation $s \to 1_{\theta_{s}} \in C^{*}(\mathcal{G}_{tight})$ is tight and any tight representation of $S$ factors through this universal one.

\begin{ppsn}
\label{main proposition}
 Let $T$ be the inverse semigroup considered in Proposition \ref{Inverse semigroup}. Denote the tight groupoid associated to $T$ by $\mathcal{G}_{tight}$. Then $\mathfrak{A}[N \rtimes H,M]$ is isomorphic to $C^{*}(\mathcal{G}_{tight})$.
\end{ppsn}
\textit{Proof.} Let $t_{a}$ and $v(m)$ be the images of $s_{a}$ and $u(m)$ in $C^{*}(\mathcal{G}_{tight})$. By Proposition \ref{tightness of the identity representation} and by the universal property of $\mathcal{G}_{tight}$, it follows that there exists a homomorphism $\rho:C^{*}(\mathcal{G}_{tight}) \to \mathfrak{A}[N \rtimes H,M]$ such that $\rho(t_{a})=s_{a}$ and $\rho(v(m))=u(m)$.

Given $a \in P$, the projections $\{u(k)e_{a}u(k)^{-1}:k \in M/M_{a}\}$ cover the projections in $T$. Since the representation of $T$ in $C^{*}(\mathcal{G}_{tight})$ is tight, it follows that 
\begin{equation*}
 \displaystyle \sum_{k \in M/M_{a}}v(k)(t_{a}t_{a}^{*})v(k)^{-1} =1
\end{equation*}

 Now the universal property of $\mathfrak{A}[N \rtimes H,M]$ implies that there exists a homomorphism $\sigma:\mathfrak{A}[N \rtimes H,M] \to C^{*}(\mathcal{G}_{tight})$ such that $\sigma(s_{a})=t_{a}$ and $\sigma(u(m))=v(m)$. It is then clear that $\sigma$ and $\rho$ are inverses of each other. This completes the proof.
\hfill $\Box$

We identify the groupoid $\mathcal{G}_{tight}$ explicitly in the rest of the article.

\section{Tight characters of the inverse semigroup $T$}
In this section, we determine the tight characters of the inverse semigroup $T$ defined in Proposition \ref{Inverse semigroup}. Let 
\begin{equation*}
  \overline{M}:=\big \{(r_{a}) \in \prod_{a \in P}M/M_{a}: r_{ab}\equiv r_{a}\mod M_{a} \big \}
\end{equation*}
We give $\overline{M}$ the subspace topology induced from the product topology on $\prod_{a \in P}M/M_{a}$. Here the finite group $M/M_{a}$ is given the discrete topology. Then $\overline{M}$ is a compact, Hausdorff topological space. Moreover $\overline{M}$ is a topological group. Note that $M$ embeds naturally into $\overline{M}$ via the imbedding $r \to (r_{a}:=r)$. The map $r \to (r_{a}:=r)$ is an imbedding since we have assumed that $\bigcap_{a \in P}M_{a}$ is trivial.

For $b \in P$ and $k \in M$, the set $U_{b,k}:=\{(r_{a}) \in \overline{M}: r_{b}\equiv k\mod M_{b} \}$ is an open set. Moreover the collection $ \{U_{b,k}: b\in P, k \in M\}$ forms a basis for $\overline{M}$. If $k \in M$ then clearly $k \in U_{b,k}$ for any $b \in P$. As a consequence, $M$ is dense in $\overline{M}$.

For $r \in \overline{M}$, let \[A_{r}:=\{f \in F: f \geq u(r_{a})e_{a}u(r_{a})^{-1} \mbox{~ for some ~} a \in P \}.\] In the next lemma, we show that for every $ r\in \overline{M}$, $A_{r}$ is an ultrafilter and all ultrafilters  are of this form.
\begin{lmma}
 \label{ultrafilter}
For $r \in \overline{M}$, $A_{r}$  is an ultrafilter. Moreover any ultrafilter is of the form $A_{r}$ for some $r \in \overline{M}$.
\end{lmma}
\textit{Proof:} Let $ r\in \overline{M}$ be given. First let us show that $A_{r}$ is a filter. Clearly $0 \notin A_{r}$. Also if $f_{1} \geq f_{2}$ and $f_{2} \in A_{r}$ then $f_{1} \in A_{r}$. Now suppose that $f_{1},f_{2} \in A_{r}$. Then there exists $a_{1},a_{2} \in P$ such that $f_{i} \geq u(r_{a_{i}})e_{a_{i}}u(r_{a_{i}})^{-1}$ for $i=1,2$. Choose $c \in a_{1}P \cap a_{2}P$. Then by Lemma \ref{decomposition}, it follows that $e_{c} \leq e_{a_{i}}$ for $i=1,2$. Since $r \in \overline{M}$, it follows that $r_{c} \equiv  r_{a_{i}} \mod M_{a_{i}}$ for $i=1,2$. Now observe that
\begin{equation*}
\begin{split}
  f_{1}f_{2} & ~\geq ~u(r_{a_{1}})e_{a_{1}}u(r_{a_{1}})^{-1}u(r_{a_{2}})e_{a_{2}}u(r_{a_{2}})^{-1} \\
             & ~= ~ u(r_{c})e_{a_{1}}u(r_{c})^{-1}u(r_{c})e_{a_{2}}u(r_{c})^{-1}\\
             & ~= ~u(r_{c})e_{a_{1}}e_{a_{2}}u(r_{c})^{-1}\\
             & ~\geq ~u(r_{c})e_{c}u(r_{c})^{-1}
\end{split}
 \end{equation*}
Thus $f_{1}f_{2} \in A_{r}$. Thus we have shown that $A_{r}$ is a filter.

Now we show $A_{r}$ is maximal. Let $A$ be a filter which contains $A_{r}$. Consider an element $f \in A$. By definition there exists $a \in P$ and scalars $\alpha_{k} \in \{0,1\}$ such that 
\[f=\sum_{k \in M/M_{a}}\alpha_{k}u(k)e_{a}u(k)^{-1}.\] But both $f$ and $u(r_{a})e_{a}u(r_{a})^{-1}$ belong to $A$ and hence their product belongs to $A$. Thus the product $fu(r_{a})e_{a}u(r_{a})^{-1}$ is non-zero. This implies that $\alpha_{r_{a}}=1$. Thus we have $f \geq u(r_{a})e_{a}u(r_{a})^{-1}$ or in other words $f \in A_{r}$. Hence $A=A_{r}$. This proves that $A_{r}$ is maximal.

Let $A$ be an ultrafilter. By Lemma \ref{maximality}, it follows that for every $a \in P$, there exists a unique $r_{a} \in M/M_{a}$ such that $u(r_{a})e_{a}u(r_{a})^{-1} \in A$. Let $r:=(r_{a})$. We claim that $r \in \overline{M}$. Let $a,b \in P$ be given. By Lemma \ref{decomposition}, we have
\begin{equation}
\label{eq}
 u(r_{a})e_{a}u(r_{a})^{-1}=\sum_{k \in M/M_{b}}u(r_{a}aka^{-1})e_{ab}u(r_{a}kak^{-1})^{-1}
\end{equation}
Since $A$ is a filter containing $u(r_{a})e_{a}u(r_{a})^{-1}$ and $u(r_{ab})e_{ab}u(r_{ab})^{-1}$, it follows that their product is non-zero. This fact together with Equation \ref{eq} implies that there exists $k \in M$, such that $r_{ab} \equiv r_{a}(aka^{-1}) \mod M_{ab}$. Thus $r_{ab}\equiv r_{a} \mod M_{a}$ for every $a,b \in P$. As a result, we have $ r \in \overline{M}$. Since $A$ is a filter it follows that $A_{r} \subset A$. We have already proved that $A_{r}$ is maximal. Thus $A=A_{r}$. This completes the proof. \hfill $\Box$ 

The following proposition identifies the tight characters of $T$.

\begin{ppsn}
\label{tight characters}
 The map $\overline{M}:r \to A_{r} \in \widehat{F_{tight}}$ is a homeomorphism.
\end{ppsn}
\textit{Proof.} It is clear from the definition that $ r \to A_{r}$ is one-one. Let us denote this map by $\phi$. We show $\phi$ is continuous. Consider a net $r^{\alpha}$  in $\overline{M}$ converging to $r$. We denote the indicator function of a set $A$ by $1_{A}$. Let $f \in F$ be given. Then there exists $a \in P$ and scalars $\alpha_{k}$ such that 
\[
 f=\sum_{k}\alpha_{k}u(k)e_{a}u(k)^{-1}
\]
Then we have
\[
 1_{A_{r^{\alpha}}}(f)=\sum_{k}\alpha_{k} \delta_{r^{\alpha}_{a},k}
\]
Since $r^{\alpha}_{a}=r_{a}$ eventually, it follows that $1_{A_{r^{\alpha}}}(f)$ converges to $1_{A_{r}}(f)$. This shows that $ r \to A_{r}$ is continuous.

Now Lemma \ref{ultrafilter} implies that $\phi$ has range $\widehat{F_{\infty}}$. Since $\overline{M}$ is compact, it follows that $\widehat{F}_{\infty}$ is compact and hence closed. Thus $\widehat{F_{\infty}}=\widehat{F_{tight}}$. Thus $\phi:\overline{M} \to \widehat{F}_{\infty}$ is one-one, onto and continuous. Since $\overline{M}$ is compact, it follows that $\phi$ is in fact a homeomorphism. This completes the proof. \hfill $\Box$

From now on we will simply denote $A_{r}$ by $r$ and $1_{A_{r}}(f)$ by $r(f)$.

\section{The groupoid $\mathcal{G}_{tight}$ of the inverse semigroup $T$}\label{dilation}

In this section, we will identify the tight groupoid $\mathcal{G}_{tight}$ associated to the inverse semigroup. Throughout this section, we assume $N=\bigcup_{a \in P}a^{-1}Ma$. By Remark \ref{key remark}, we can very well assume this. There is another natural groupoid which arises out of the following construction.

For every $a \in P$, the co-isometry $s_{a}^{*}$ will give rise to an injection on $\overline{M}$ and the unitary $u(m)$ for $m \in M$ will act as a bijection  on $\overline{M}$. Thus we get an action of the semigroup $M \rtimes P$, as injections, on $\overline{M}$. Now the space $\overline{M}$ can be enlarged to a space $\overline{N}$ and the action of $M \rtimes P$ can be dilated to get an action of $G= N \rtimes H$ on $\overline{N}$. We can then consider the transformation groupoid $\overline{N} \rtimes G$. But the unit space of $\mathcal{G}_{tight}$ is $\overline{M}$. Thus we restrict the transformation groupoid $\overline{N} \rtimes G$ to $\overline{M}$ and prove that it is isomorphic to $\mathcal{G}_{tight}$. 

This dilation procedure  has appeared in several works  [See \cite{Laca-dilation}, \cite{Larsen-Li} ]. The basic principle goes back to \cite{Ore}. 

  First let us explain the action of $M \rtimes P$ on $\overline{M}$. The action of $M$ on $\overline{M}$ is by left multiplication as $M$ is a subgroup of $\overline{M}$. Let $a \in P$  and $r \in \overline{M}$ be given. For $b \in P$, choose $c \in aP \cap bP$ and write $c$ as $c=a\alpha=b \beta$. We will use the notation as in Remark \ref{conjugation by s}. Note that $M_{c} \subset M_{b}$ and we denote the induced quotient map $M/M_{c} \to M/M_{b}$ by $q_{b,c}$. Define $m_{b}=q_{b,c}(\pi_{\alpha}^{a}(r_{\alpha}))$. First let us show that $m_{b}$ depends only on $a$ and $b$ and not on the choices made. 

Suppose $c_{1}=a\alpha_{1}=b\beta_{1}$ and $c_{2}=a\alpha_{2}=b\beta_{2}$. Choose $\gamma_{1},\gamma_{2} \in P$ such that $\alpha_{1}\gamma_{1}=\alpha_{2}\gamma_{2}$. Note that this implies $c_{1}\gamma_{1}=c_{2}\gamma_{2}$. 
Now we have
\begin{equation*}
 \begin{split}
 q_{b,c_{i}}\pi_{\alpha_{i}}^{a}(r_{\alpha_{i}})& = q_{b,c_{i}}\Big(\pi_{\alpha_{i}}^{a}\big(q_{\alpha_{i},\alpha_{i}\gamma_{i}}(r_{\alpha_{i}\gamma_{i}})\big)\Big) \\ 
                                                &=q_{b,c_{i}}\Big(q_{c_{i},c_{i}\gamma_{i}}\big(\pi_{\alpha_{i}\gamma_{i}}^{a}(r_{\alpha_{i}\gamma_{i}})\big)\Big)\\
                                                &=q_{b,c_{i}\gamma_{i}}\big(\pi_{\alpha_{i}\gamma_{i}}^{a}(r_{\alpha_{i}\gamma_{i}})\big).                                                
 \end{split}
\end{equation*}
Note that the right hand side is constant for $i=1,2$. Thus we have \[q_{b,c_{1}}(\pi_{\alpha_{1}}^{a}(r_{\alpha_{1}}))=q_{b,c_{2}}(\pi_{\alpha_{2}}^{a}(r_{\alpha_{2}})).\]
This shows that $m_{b}$ is well defined. We leave it to the reader to check that $\tilde{m}=(m_{b}) \in \overline{M}$.

 On $M$, the action of $P$ is the usual conjugation. From now on, we denote the element $\tilde{m}$ by $ara^{-1}$. This way $P$ acts on $\overline{M}$ injectively and continuously. This action of $P$ together with the left multiplication action of $M$ defines an action of $M \rtimes P$ on $\overline{M}$ (as injective,continuous transformations). We leave the details to the reader.

\begin{lmma}
\label{kernel}
For $a \in P$, the kernel of the projection map $\overline{M} \ni (y_{b}) \to y_{a} \in M/M_{a}$ is $a\overline{M}a^{-1}$.
\end{lmma}
\textit{Proof.} By definition, it follows that $a\overline{M}a^{-1}$ is in the kernel of the $a^{\mbox{th}}$ projection. Now let $y=(y_{b})$ be such that $y_{a}=1$. Since $M$ is dense in $\overline{M}$, there exists a sequence $y^{n} \in M$ such that $y^{n} \to y$ in $\overline{M}$. As $M/M_{a}$ is finite, we can without loss of generality assume that $y^{n} \in M_{a}$ for every $n$. Thus there exists $x^{n} \in M$ such that $y^{n}=ax^{n}a^{-1}$. But $\overline{M}$ is compact. Thus, by passing to a subsequence if necessary, we can assume that $x^{n}$ converges to an element say $x \in \overline{M}$. Since conjugation by $a$ is continuous, it follows that $y^{n}=ax^{n}a^{-1}$ converges to $axa^{-1}$. But $y^{n}$ converges to $y$. Thus $axa^{-1}=y$. This completes the proof. \hfill $\Box$

Now let us explain the dilation procedure that we promised at the beginning of this section.
Consider the set $\overline{M} \times P$ and define a relation on $\overline{M} \times P$ by $(x,a) \sim (y,b)$ if there exists $\alpha,\beta \in P$ such that 
$\alpha a=\beta b$ and $\alpha x \alpha^{-1}=\beta y \beta^{-1}$. We leave the following routine checking to the reader.
\begin{enumerate}
 \item The relation $\sim$ is an equivalence relation. We denote the equivalence class containing $(x,a)$ by $[(x,a)]$.
 \item Let $\overline{N}:=\overline{M} \times P/\sim$. Then $\overline{N}$ is a  group. The multiplication on $\overline{N}$ is defined as follows. For $a,b \in P$, choose $\alpha$ and $\beta$ such that $\alpha a=\beta b$.
Then \[
      [(x,a)][(y,b)]=[(\alpha x \alpha^{-1}\beta y \beta^{-1},\alpha a)]
     \]
The identity element of $\overline{N}$ is $[(e,e)]$ where $(e,e)$ is the identity element of $\overline{M}\times P$ and the inverse of $[(x,a)]$ is $[(x^{-1},a)]$.
\item The group $\overline{N}$ is a locally compact Hausdorff topological group when $\overline{N}$ is given the quotient topology. Here  $P$ is given the discrete topology. 
 \item The map $M\ni x \to [(x,e)] \in \overline{N}$ is a topological embedding. Thus $\overline{M}$ can be viewed as a subset of $\overline{N}$. Moreover $\overline{M}$ is a compact open subgroup of $\overline{N}$.
\item The map $N \ni a^{-1}ma \to [(m,a)] \in \overline{N}$ is an embedding. When $N$ is viewed as a subset of $\overline{N}$ via this embedding, $N$ is dense in $\overline{N}$. Also $N \cap \overline{M}=M$.
\item Let $a \in P$ be given. Define a map $\phi_{a}:\overline{N} \to \overline{N}$ as follows. Given $[(x,b)] \in \overline{N}$, choose $\alpha,\beta \in P$ such that $\alpha a=\beta b$. Define $\phi_{a}([(x,b)])=[\beta x \beta^{-1},\alpha)]$. One checks that $\phi_{a}$ is well defined. Moreover for $a \in P$, $\phi_{a}$ is a homeomorphism with $\phi_{a}^{-1}$ given by $\phi_{a}^{-1}[(x,b)]=[(x,ba)]$. Note that $\phi_{a}$ restricted to $N$ is the usual conjugation. Also $\phi_{a}\phi_{b}=\phi_{ab}$ for $a,b \in P$. For $m \in M$ , define $\psi_{m}:\overline{N} \to \overline{N}$ as $\psi_{m}([(x,a)])=[(ama^{-1}x,a)]$. That is $\psi_{m}$ is just left multiplication by $m$. One also has the following commutation relation. For $a \in P$ and $m \in M$, 
\[
\phi_{a}\psi_{m}=\psi_{ama^{-1}}\phi_{a}.
\]
\item Since we have assumed that $N=\bigcup_{a \in P}a^{-1}Ma$, it follows that any element of $g \in G=N \rtimes H$ can be written as $g=a^{-1}mb$ with $a,b \in P$ and $m \in M$. The map $a^{-1}mb \to \phi_{a}^{-1}\psi_{m}\phi_{b}$ is well defined and defines an action of $G$ on $\overline{N}$. If $h=a^{-1}b \in H$ and $x \in \overline{N}$, we denote $\phi_{a}^{-1}\phi_{b}(x)$ as $hxh^{-1}$. If $n=a^{-1}ma$ and $x\in \overline{N}$, we denote $\phi_{a}^{-1}\psi_{m}\phi_{a}(x)$ as $nx$.
\item Note that $\overline{N}=\bigcup_{a \in P}a^{-1}\overline{M}a$.
\item \label{universal property}  \textbf{Universal Property:} Let $L$ be a locally compact Hausdorff topological group on which $H$ acts by group homomorphism. Suppose that $K$ is a compact open subgroup of $L$ which is invariant under $P$ and $L=\bigcup_{a \in P}a^{-1}K$. If $\phi:\overline{M} \to K$ is a $P$-equivariant continuous bijection then the map $\overline{N} \ni a^{-1}xa \to a^{-1}.\phi(x) \in L$ is a topological isomorphism and is $H$ -equivariant.
\end{enumerate}

\begin{rmrk}
 It is not difficult to show by using (9) that $\overline{N}$ is the pro-finite completion of $N$ when $N$ is given the topology induced by the neighbourhood base $\{aMa^{-1}:a \in H\}$ at the identity. In \cite{Quigg-Landstad}, the pro-finite completion model of $\overline{N}$ is used.
\end{rmrk}

 When considering transformation groupoids, we consider only right actions of groups and thus we change the above left action of $G$ on $\overline{N}$ to a right action simply by defining $x.g=g^{-1}x$ for $x\in \overline{N}$ and $g \in G$. Now consider the transformation groupoid $\overline{N} \rtimes G$ and restrict it to $\overline{M}$. We show that the groupoid $\mathcal{G}_{tight}$ of the inverse semigroup $T$ is isomorphic to the groupoid $\overline{N} \rtimes G|_{\overline{M}}$ i.e. to the transformation groupoid $\overline{N} \rtimes G$ restricted to the unit space $\overline{M}$. We will start with two lemmas which will be extremely useful to prove this.

\begin{lmma}
\label{Universal group of the inverse semigroup}
If $a_{1}^{-1}m_{1}b_{1}=a_{2}^{-1}m_{2}b_{2}$ then $s_{a_{1}}^{*}u(m_{1})s_{b_{1}}=s_{a_{2}}^{*}u(m_{2})s_{b_{2}}$. 
\end{lmma}
\textit{Proof.} Suppose $a_{1}^{-1}m_{1}b_{1}=a_{2}^{-1}m_{2}b_{2}$. Then $a_{1}^{-1}m_{1}a_{1}=a_{2}^{-1}m_{2}a_{2}$ and $a_{1}^{-1}b_{1}=a_{2}^{-1}b_{2}$. Choose $\beta_{1},\beta_{2} \in P$ such that $\beta_{1}b_{1}=\beta_{2}b_{2}$. Then  $a_{1}a_{2}^{-1}=\beta_{1}^{-1}\beta_{2}=b_{1}b_{2}^{-1}$. Hence $\beta_{1}m_{1}\beta_{1}^{-1}=\beta_{2}m_{2}\beta_{2}^{-1}$. Now observe that 
\begin{equation*}
\begin{split}
s_{a_{1}}^{*}u(m_{1})s_{b_{1}}&=s_{a_{1}}^{*}u(m_{1})s_{\beta_{1}}^{*}s_{\beta_{1}}s_{b_{1}} \\
                                                       &=s_{a_{1}}^{*}s_{\beta_{1}}^{*}u(\beta_{1}m_{1}\beta_{1}^{-1})s_{\beta_{1}b_{1}}\\
                                                      &=s_{\beta_{1}a_{1}}^{*}u(\beta_{1}m_{1}\beta_{1}^{-1})s_{\beta_{1}b_{1}} \\
                                                      &=s_{\beta_{2}a_{2}}^{*}u(\beta_{2}m_{2}\beta_{2}^{-1})s_{\beta_{2}b_{2}}\\
                                                      &=s_{a_{2}}^{*}s_{\beta_{2}}^{*}u(\beta_{2}m_{2}\beta_{2}^{-1})s_{\beta_{2}b_{2}}\\
                                                     &=s_{a_{2}}^{*}u(m_{2})s_{\beta_{2}}^{*}s_{\beta_{2}}s_{b_{2}} \\
                                                     &=s_{a_{2}}^{*}u(m_{2})s_{b_{2}} 
\end{split}
\end{equation*} 
This completes the proof. \hfill $\Box$

\begin{lmma}
 \label{omitting projections}
In $\mathcal{G}_{tight}$, $[(r,s_{a}^{*}u(m)fu(n)s_{b})]=[(r,s_{a}^{*}u(mn)s_{b})]$.
\end{lmma}
\textit{Proof.} First observe that $[(r,s_{a}^{*})][r.s_{a}^{*},u(m)fu(n)s_{b}]=[(r,s_{a}^{*}u(m)fu(n)s_{b})]$. Thus it is enough to consider the case when $a$ is the identity element of $P$. Now let $s=u(m)fu(n)s_{b}$, $t=u(mn)s_{b}$ and $e=u(m)fu(m)^{-1}$. Observe that $s=et$. Thus $ss^{*}=ett^{*}e$. Hence $r(ss^{*})=1$ implies $r(e)=1$ and $r(tt^{*})=1$. Moreover $es=s=et$. Thus $[(r,s)]=[(r,t)]$. This completes the proof. \hfill $\Box$

Now we can state our main theorem.

\begin{thm}
 \label{main theorem}
Let $\phi:\overline{N} \rtimes G|_{\overline{M}} \to \mathcal{G}_{tight}$ be the map defined by 
\[
 \phi\big((x,a^{-1}mb)\big)=[(x,s_{a}^{*}u(m)s_{b})].
\]
Then $\phi$ is a topological groupoid isomorphism.
\end{thm}
\textit{Proof.} First let us show that $\phi$ is well defined. Let $(x,a^{-1}mb) \in \overline{N} \rtimes G|_{\overline{M}}$. Then by definition, there exists $y \in \overline{M}$ such that $m^{-1}axa^{-1}=byb^{-1}$. Choose $\alpha$ and $\beta$ in $P$ such that $c:=a\alpha=b\beta$.  By definition, this means that 
$\pi_{\alpha}^{a}(x_{\alpha})\equiv q_{c}(m)\pi_{\beta}^{b}(y_{\beta})$. Now Remark \ref{conjugation by s} implies that 
\[s_{a}^{*}u(m)e_{b}u(m)^{-1}s_{a} \geq u(x_{\alpha})e_{\alpha}u(x_{\alpha})^{-1}.\] Hence $x(s_{a}^{*}u(m)e_{b}u(m)^{-1}s_{a})=1$. Thus we have shown that $\phi$ is well-defined.

Before we show $\phi$ is a surjection,  let us show that if $[(x,s_{a}^{*}u(m)s_{b})] \in \mathcal{G}_{tight}$ then $(x,a^{-1}mb) \in \overline{N} \rtimes G|_{\overline{M}}$. To that effect, assume that $x(s_{a}^{*}u(m)e_{b}u(m)^{-1}s_{a})=1$. Choose $c \in aP\cap bP$ and write $c=a\alpha=b\beta$. By Remark \ref{conjugation by s}, it follows that there exists $y \in M/M_{\beta}$ such that $q_{c}(m^{-1})\pi_{\alpha}^{a}(x_{\alpha})=\pi_{\beta}^{b}(y)$. This implies that the $b^{th}$ co-ordinate of $m^{-1}axa^{-1}$ is $1$ i.e. the identity element of $M/M_{b}$. Now Lemma \ref{kernel} implies that there exists $z \in \overline{M}$ such that $m^{-1}axa^{-1}=bzb^{-1}$. Hence $(x,a^{-1}mb) \in \overline{N} \rtimes G|_{\overline{M}}$. Surjectivity is then an immediate consequence of Lemma  \ref{omitting projections}.

Now we show $\phi$ is injective. Suppose $[(x,s_{a_{1}}^{*}u(m_{1})s_{b_{1}})]=[(x,s_{a_{2}}^{*}u(m_{2})s_{b_{2}})]$. Then there exists a projection $e \in F$  such that $0 \neq e(s_{a_{1}}^{*}u(m_{1})s_{b_{1}})=e(s_{a_{2}}^{*}u(m_{2})s_{b_{2}})$. We can without loss of generality assume that $e=u(r_{c})e_{c}u(r_{c})^{-1}$. By Remark \ref{good remark} and by reading the above equality in the regular representation, we immediately obtain $a_{1}^{-1}b_{1}=a_{2}^{-1}b_{2}$ and $a_{1}^{-1}m_{1}b_{1}=a_{2}^{-1}m_{2}b_{2}$. This implies that $\phi$ is injective.

Now let us show that $\phi$ is a groupoid morphism. First we show that $\phi$ preserves the range and source. By definition, $\phi$ preserves the range.  Observe that $\phi$ is continuous and this is a direct consequence of Proposition \ref{tight characters}. Let $\gamma=(x,a^{-1}mb) \in \overline{N} \rtimes G|_{\overline{M}}$. Since $M$ is dense in $\overline{M}$ there exists a sequence $x_{n} \in M$ such that $x_{n}$ converges to $x$. Moreover the action of $G$ on $\overline{N}$ is continuous and $\overline{M}$ is compact and open. Thus we can assume that $(x_{n},a^{-1}mb) \in \overline{N} \rtimes G|_{\overline{M}}$ for every $n$. By definition, there exists $y \in \overline{M}$ such that $axa^{-1}=mbyb^{-1}$. Also let $y_{n}$ be such that $ax_{n}a^{-1}=mby_{n}b^{-1}$. 

To keep things clear, if $z \in \overline{M}$, we denote the character determined by $z$ as $\xi_{z}$. Let $v:=s_{a}^{*}u(m)s_{b}$.
Now if can show that $\xi_{x_{n}}.v=\xi_{y_{n}}$ then it will follow from continuity of $\phi$ that $\xi_{x}.v=\xi_{y}$. Thus we only need to show that $s(\phi(\gamma))=\phi(s(\gamma))$ for $\gamma=(x,a^{-1}mb)$ with $x\in M$.

Now let $(x,a^{-1}mb) \in \overline{N} \rtimes G|_{\overline{M}}$ with $x \in M$. Then there exists $y \in M$ such that 
$axa^{-1}=mbyb^{-1}$. Let $v=s_{a}^{*}u(m)s_{b}$.  To show $\xi_{x}.v=\xi_{y}$, as $\xi_{y}$ is maximal, it is enough to show that the support of $\xi_{y}$ is contained in $\xi_{x}.v$. Again it is enough to show that $u(y)e_{c}u(y)^{-1}$ is in the support of $\xi_{x}.v$. Choose $\alpha,\beta$ such that $a\alpha=bc\beta$.  Note that
\begin{equation*}
\begin{split}
vu(y)e_{c}u(y)^{-1}v^{*}&=s_{a}^{*}u(m)s_{b}u(y)e_{c}u(y)^{-1}s_{b}^{*}u(m)^{-1}s_{a}\\
                           &=s_{a}^{*}u(mbyb^{-1})s_{b}e_{c}s_{b}^{*}u(mbyb^{-1})^{-1}s_{a}\\
                           &=s_{a}^{*}u(axa^{-1})e_{bc}u(axa^{-1})^{-1}s_{a} \\
                          &=u(x)s_{a}^{*}e_{bc}s_{a}u(x)^{-1} \\
                         & \geq u(x)s_{a}^{*}e_{bc\beta}s_{a}u(x)^{-1}\\
                         &=u(x)s_{a}^{*}e_{a\alpha}s_{a}u(x)^{-1} \\
                         &=u(x)e_{\alpha}u(x)^{-1} \in \mbox{supp}(\xi_{x})
\end{split}
  \end{equation*}
Hence $u(y)e_{c}u(y)^{-1}$ is in the support of $\xi_{x}.v$. Thus we have shown that $\xi_{x}.v=\xi_{y}$. This proves that $\phi$ preserves the source.

Now we show $\phi$ preserves multiplication. Let $\gamma_{1}=(x_{1},a_{1}^{-1}m_{1}b_{1})$ and $\gamma_{2}=(x_{2},a_{2}^{-1}m_{2}b_{2})$. Since $\phi$  preserves the range and source, it follows that $\gamma_{1}$ and $\gamma_{2}$ are composable if and only if $\phi(\gamma_{1})$ and $\phi(\gamma_{2})$ are composable. Choose $\alpha,\beta \in P$ such that $\beta b_{1}=\alpha a_{2}$.  . Now 
\begin{align*}
\phi(\gamma_{1})\phi(\gamma_{2})&=[(x_{1},s_{a_{1}}^{*}u(m_{1})s_{b_{1}}s_{a_{2}}^{*}u(m_{2})s_{b_{2}})]\\
                                                         &=[(x_{1},s_{\beta a_{1}}^{*}u(\beta m_{1}\beta^{-1})e_{\alpha a_{2}}u(\alpha m_{2} \alpha^{-1})s_{\alpha b_{2}})]  ~\big(\mbox { by Eq. \ref{commutation rule}}\Big) \\
                                                         &=[(x_{1},s_{\beta a_{1}}^{*}u(\beta m_{1} \beta^{-1} \alpha m_{2} \alpha^{-1})s_{\alpha b_{2}})] ~\big( \mbox{ by Remark \ref{omitting projections}} \Big)\\
                                                         &=\phi(\gamma_{1}\gamma_{2}).
\end{align*}
It is easily verifiable that $\phi$ preserves inversion.

 For an open subset $U$ of $\overline{M}$ and $g=a^{-1}mb$, consider the open set 
\[
 \theta(U,g):=\{x \in \overline{M}:x.g \in \overline{M}\}
\]
The collection $\{\theta(U,g)\}$ forms a basis for $\overline{N} \rtimes G|_{\overline{M}}$. Moreover $\phi(\theta(U,g))=\theta(U,s_{a}^{*}u(m)s_{b})$. Thus $\phi$ is an open map. Thus we have shown that $\phi$ is a homeomorphism. This completes the proof. \hfill $\Box$

\begin{crlre}
\label{main corollary}
 The algebra $\mathfrak{A}[N \rtimes H,M]$ is isomorphic to $C^{*}(\overline{N} \rtimes G|_{\overline{M}})$. 
\end{crlre}
\textit{Proof.} This follows from Theorem \ref{main theorem} and Proposition \ref{main proposition}. \hfill $\Box$

\section{Simplicity of $\mathfrak{A}_{r}[N \rtimes H,M]$}
Let us recall a few definitions from \cite{Claire}. Let $\mathcal{G}$ be an r-discrete groupoid and we denote its unit space by $\mathcal{G}^{0}$. The relation $\sim$ defined by $x \sim y$ if and only if there exists $\gamma \in \mathcal{G}$ such that $s(\gamma)=x$ and $r(\gamma)=y$  is an equivalence relation on $\mathcal{G}^{0}$. A subset $E \subset \mathcal{G}^{0}$ is said to be invariant if given $x \in E$ and $y \sim x$ then $y \in E$. For $x \in \mathcal{G}$, let $\mathcal{G}(x):=\{\gamma \in \mathcal{G}:s(\gamma)=r(\gamma)=x\}$ be the isotropy group of $x$.

A subset $S \subset \mathcal{G}$ is said to be a bi-section if the range and source maps restricted to $S$ are one-one. If $S$ is a bisection, let $\alpha_{S}:r(S)\to s(S)$ be defined by $\alpha_{S}:=s \circ r^{-1}$. 

 The groupoid $\mathcal{G}$ is said to be 

\begin{itemize}
\item  minimal if the only non-empty, open invariant subset of $\mathcal{G}^{0}$ is $\mathcal{G}^{0}$.
\item topologically principal if the set of $x \in \mathcal{G}^{0}$ for which $\mathcal{G}(x)=\{x\}$ is dense in $\mathcal{G}^{0}$.
\item locally contractive if  for every non-empty open subset $U$ of $\mathcal{G}^{0}$, there exists an open subset $V \subset U$ and an open bisection $S$ with $\overline{V} \subset s(S)$ and $\alpha_{S^{-1}}(\overline{V})$ not contained in $V$.
\end{itemize}

Conjugation by $P$ on $M$ gives rise to a semigroup homomorphism from $P$ to the semigroup of injective maps on $M$. In \cite{Quigg-Landstad}, the action of $P$ on $M$ is called an effective action if the above semigroup homomorphism is injective i.e. given $h \in H$ with $h\neq 1$, then there exists $s \in M$ such that $hsh^{-1}\neq s$. In \cite{Quigg-Landstad}, the following facts were proved about the transformation groupoid $\overline{N} \rtimes G$. 
\begin{enumerate}
\item The groupoid  $\overline{N} \rtimes G$ is minimal and locally contractive.
\item The groupoid  $\overline{N} \rtimes G$ is topologically principal if and only if $P$ acts effectively on $M$.
\item Thus the reduced $C^{*}$-algebra $C_{red}^{*}(\overline{N} \rtimes G)$ is simple and purely infinite if $P$ acts effectively on $M$. [Refer to \cite{Claire}].
\end{enumerate} 
Analogous statements hold for the groupoid $\mathcal{G}_{tight}$ associated to the inverse semigroup $T$. 
\begin{rmrk}
In \cite{Quigg-Landstad}, only the if part (in (2)) was proved. But then the other direction i.e. if $\overline{N} \rtimes G$ is topologically principal then $P$ acts effectively on $M$ is easy to verify.
\end{rmrk}

Also note that $\overline{M}$ is a closed  subset of $\overline{N}$ which meets each $G$ orbit of $\overline{N}$. Moreover $\overline{M}$ is open as well. Hence by appealing to Example 2.7 in \cite{MRW} , we conclude that $C^{*}(\overline{N} \rtimes G)$ and $C^{*}(\overline{N} \rtimes G|_{\overline{M}})$ are Morita-equivalent.

We end this section by showing that $\mathfrak{A}_{r}[N \rtimes H,M]$ is isomorphic to the reduced $C^{*}$-algebra $C_{red}^{*}(\mathcal{G}_{tight})$.

\begin{ppsn}
\label{reduced}
Let $\mathcal{G}:=\overline{N} \rtimes G|_{\overline{M}}$. Then the reduced $C^{*}$-algebra of the groupoid $\mathcal{G}$ is isomorphic to $\mathfrak{A}_{r}[N \rtimes H,M]$.
\end{ppsn}
\textit{Proof.} Let $e$ be the identity element of $\overline{M}$. Define $\mathcal{G}^{e}:=\{\gamma \in \mathcal{G}:r(\gamma)=e\}$. Then $\mathcal{G}^{e}:=\{(e,hm): m \in M,h \in H\}$. Thus $L^{2}(\mathcal{G}^{e})$ can be identified with $\ell^{2}(M)\otimes \ell^{2}(H)$. Consider the representation $\pi_{e}$ of $C_{red}^{*}(\mathcal{G})$ on $L^{2}(\mathcal{G}^{e})$ defined as follows. For $f \in C_{c}(\mathcal{G})$, define $\pi_{e}(f)$ by the following formula.
\[
(\pi_{e}(f)(\xi))(\gamma):=\displaystyle \sum_{\gamma_{1} \in \mathcal{G}^{e}}f(\gamma^{-1}\gamma_{1})\xi(\gamma_{1})
\]
Since $M$ is dense in $\overline{M}$, it follows that the largest open invariant set not containing $e$ is the empty set. Hence $\pi_{e}$ is faithful. 

For $a \in P$ and $m \in M$, we let $S_{a}$ and $U(m)$ be the images of $s_{a}$ and $u(m)$ in $C_{red}^{*}(\mathcal{G})$. Let $\{\delta_{m}\otimes \delta_{b}: m \in M, b \in H\}$ be the canonical basis of $\ell^{2}(M)\otimes \ell^{2}(H)$. Consider the unitary operator $V$ on $\ell^{2}(M)\otimes \ell^{2}(H)$ defined by 
\[
V(\delta_{m}\otimes \delta_{b}):=\delta_{m^{-1}}\otimes \delta_{b^{-1}}
\]
For $a \in P$ and $k \in M$, we leave it to the reader to check the following equality.
\begin{equation*}
\begin{split}
V\pi_{e}(S_{a})V^{*}(\delta_{m}\otimes \delta_{b})&=\delta_{ama^{-1}}\otimes \delta_{ab} \\
V\pi_{e}(U(k))V^{*}(\delta_{m}\otimes \delta_{b})&=\delta_{km}\otimes \delta_{b}
\end{split}
\end{equation*}
Since $\{S_{a}:a \in P\}$ and $\{U(k):k \in M\}$ generate  $C_{red}^{*}(\mathcal{G})$, it follows that $C_{red}^{*}(\mathcal{G})$  is isomorphic to $\mathfrak{A}_{r}[N \rtimes H,M]$. This completes the proof. \hfill $\Box$

We now show that Corollary \ref{main corollary} and Proposition \ref{reduced} can also be expressed in terms of  crossed products as in \cite{Quigg-Landstad}. We need to digress a bit before we do this. 

Let $\mathcal{G}$ be an $r$-discrete, locally compact and Hausdorff groupoid.
Let $Y \subset \mathcal{G}^{0}$ be a compact open subset of the unit space. Assume that $Y$ meets each orbit of $\mathcal{G}^{0}$. Let 
\begin{align*}
 \mathcal{G}^{Y}:&=\{\gamma \in \mathcal{G}:s(\gamma) \in Y \}\\
 \mathcal{G}^{Y}_{Y}:&=\{\gamma \in \mathcal{G}:s(\gamma),r(\gamma) \in Y \}
\end{align*}
Since $Y$ is clopen, it follows that $\mathcal{G}^{Y}$ and $\mathcal{G}_{Y}^{Y}$ are clopen. Thus if $f \in C_{c}(G^{Y})$, then $f$ can be extended to an element in $C_{c}(\mathcal{G})$ by declaring its value to be zero outside $\mathcal{G}^{Y}$. Thus we have the inclusion $C_{c}(\mathcal{G}^{Y}) \subset C_{c}(\mathcal{G})$. Similarly, we have the inclusion $C_{c}(\mathcal{G}^{Y}_{Y}) \subset C_{c}(\mathcal{G}^{Y})$. The algebra $C_{c}(\mathcal{G}^{Y}_{Y})$ is a $*$-subalgebra of $C_{c}(\mathcal{G})$. 

The space $C_{c}(\mathcal{G}^{Y})$ is a pre-Hilbert $C_{c}(\mathcal{G}^{Y}_{Y}) \subset C^{*}(\mathcal{G}^{Y}_{Y})$ module with the inner product and the right multiplication given by
\begin{align*}
 <f_{1},f_{2}>(\gamma)&=\sum_{\gamma_{1}\gamma_{2}=\gamma}\overline{f_{1}(\gamma_{1}^{-1})}f_{2}(\gamma_{2}) ~ \text{~~~for~~}\gamma \in \mathcal{G}^{Y}_{Y},~ f_{1},f_{2} \in C_{c}(\mathcal{G}^{Y})\\
 (f.g)(\gamma)&=\sum_{\gamma_{1}\gamma_{2}=\gamma}f(\gamma_{1})g(\gamma_{2})  \text{~~~for~~} \gamma \in \mathcal{G}^{Y}, ~f \in C_{c}(\mathcal{G}^{Y}),~g \in C_{c}(\mathcal{G}_{Y}^{Y})
\end{align*}
Moreover there is left action of $C_{c}(\mathcal{G})$ on $C_{c}(\mathcal{G}^{Y})$ and it is  given by 
\begin{align*}
 (f.\phi)(\gamma)&= (f*\phi)(\gamma) \\
                 &=\sum_{\gamma_{1}\gamma_{2}=\gamma}f(\gamma_{1})\phi(\gamma_{2})
\end{align*}
for $\gamma \in \mathcal{G}^{Y}$, $f \in C_{c}(\mathcal{G})$ and $\phi \in C_{c}(\mathcal{G}^{Y})$. 

Now Theorem 2.8 and Example 2.7 of \cite{MRW} implies the following. The ``completion'' of $C_{c}(\mathcal{G})$-$C_{c}(\mathcal{G}^{Y}_{Y})$ bimodule $C_{c}(\mathcal{G}^{Y})$ is a $C^{*}(\mathcal{G})$-$C^{*}(\mathcal{G}^{Y}_{Y})$ imprimitivity bimodule implementing a strong Morita equivalence between $C^{*}(\mathcal{G})$ and $C^{*}(\mathcal{G}^{Y}_{Y})$.

Let us denote the completion of $C_{c}(\mathcal{G}^{Y})$ by $\mathcal{E}$.  For $x,y \in \mathcal{E}$, let $\theta_{x,y}$ be the compact operator on $\mathcal{E}$ defined by $\theta_{x,y}(z)=x<y,z>$. For $x \in \mathcal{E}$, the operator norm of $\theta_{x,x}$ is $||x||^{2}$.

The following proposition  has also appeared in  \cite{Li-nuclearity}. (See Lemma 5.18 in \cite{Li-nuclearity}.) The proof is exactly as in \cite{Li-nuclearity}. We include the proof for the sake of completeness.

\begin{ppsn}
 \label{Isometric embedding}
The inclusion $C_{c}(\mathcal{G}^{Y}_{Y}) \subset C_{c}(\mathcal{G})$ extends to an isometric embedding from $C^{*}(\mathcal{G}^{Y}_{Y})$ to $C^{*}(\mathcal{G})$. Also the inclusion $C_{c}(\mathcal{G}^{Y}_{Y}) \subset C_{c}(\mathcal{G})$ extends to an isometric embedding from $C^{*}_{red}(\mathcal{G}_{Y}^{Y})$ to $C^{*}_{red}(\mathcal{G})$.
\end{ppsn}
\textit{Proof.} Let $f \in C_{c}(\mathcal{G}^{Y}_{Y})$ be given. Consider $f$ as an element of $C_{c}(\mathcal{G}^{Y}) \subset \mathcal{E}$. Then $\theta_{f,f}$ restricted to $C_{c}(\mathcal{G}^{Y})$ is just multiplication by $f*f^{*}$. Since $\mathcal{E}$ is a $C^{*}(\mathcal{G})$-$C^{*}(\mathcal{G}^{Y}_{Y})$ imprimitivity bimodule, it follows that  
\begin{align*}
||f||^{2}_{C^{*}(\mathcal{G})}&=||f * f^{*}||_{C^{*}(\mathcal{G})}\\
                              &=||\theta_{f,f}|| \\
                              &=||f||^{2}_{\mathcal{E}} \\
                              &=||f^{*}*f||_{C^{*}(\mathcal{G}^{Y}_{Y})} \\
                              &=||f||^{2}_{C^{*}(\mathcal{G}^{Y}_{Y})}
                              \end{align*}

For $x \in \mathcal{G}^{0}$, let $\mathcal{G}^{(x)}:=r^{-1}(x)$. Consider $\ell^{2}(\mathcal{G}^{(x)})$ and let $\{\delta_{\gamma}:\gamma \in \mathcal{G}^{(x)}\}$ be the standard orthonormal basis. Consider the representation $\pi_{x}$ of $C_{c}(\mathcal{G})$ on $\ell^{2}(\mathcal{G}^{(x)})$ defined by 
\begin{equation}
\label{reduced representation}
 \pi_{x}(f)(\delta_{\gamma})=\sum_{\alpha \in \mathcal{G}^{(x)}}f(\alpha^{-1}\gamma)\delta_{\alpha}.
\end{equation}
The reduced $C^{*}$-algebra $C^{*}_{red}(\mathcal{G})$ is the completion of $C_{c}(\mathcal{G})$ under the norm $||.||$ given by $|f||_{red}=sup_{x \in \mathcal{G}^{0}}||\pi_{x}(f)||$. (We refer the reader to \cite{Ren_book}.)

Let $\mathcal{G}^{(x)}_{Y}:=\{\gamma \in \mathcal{G}^{(x)}: s(\gamma) \in Y \}$. If $x \in Y$, let $\pi_{x}^{Y}$ be the representation of $C_{c}(\mathcal{G}^{Y}_{Y})$ on $\ell^{2}(\mathcal{G}^{(x)}_{Y})$ defined by the same formula as in Eq. \ref{reduced representation}. Now observe the following.

\begin{enumerate}
 \item  Let $\gamma_{0} \in \mathcal{G}$ be such that $s(\gamma_{0})=x$ and $r(\gamma_{0})=y$. Then $U:\ell^{2}(\mathcal{G}^{(x)}) \to \ell^{2}(\mathcal{G}^{(y)})$ defined by $\displaystyle U(\delta_{\gamma})=\delta_{\gamma_{0}\gamma}$ is a unitary. Moreover $U\pi_{x}(.)U^{*}=\pi_{y}(.)$.
\item Since $Y$ meets each orbit of $\mathcal{G}^{0}$, it follows from $(1)$ that for $f \in C_{c}(\mathcal{G})$, $\displaystyle ||f||_{red}=\sup_{x \in Y} ||\pi_{x}(f)||$.
\item If $x \in Y$, then write $\ell^{2}(\mathcal{G}^{(x)})$ as $\ell^{2}(\mathcal{G}^{(x)})=\ell^{2}(\mathcal{G}^{(x)}_{Y})\oplus (\ell^{2}(\mathcal{G}^{(x)}_{Y}))^{\perp}$. With this decomposition, for $f \in C_{c}(\mathcal{G}^{Y}_{Y})$, we have $\pi_{x}(f)=\pi_{x}^{Y}(f) \oplus 0$.
\end{enumerate}

Now the above three observations imply that for $f \in C_{c}(\mathcal{G}^{Y}_{Y})$, $||f||_{C^{*}_{red}(\mathcal{G}^{Y}_{Y})}=||f||_{C^{*}_{red}(\mathcal{G})}$. This completes the proof. \hfill $\Box$

\begin{rmrk}
 The representations used to define the regular representation in \cite{Ren_book} is different from what we have used. But the inversion map of the groupoid intertwines our representations with those used in \cite{Ren_book}. 
\end{rmrk}

The $C^{*}$-algebra of the groupoid $\overline{N} \rtimes G$ is naturally isomorphic to $C_{0}(\overline{N}) \rtimes G$ . Let $\varPhi:C_{c}(\overline{N}) \rtimes G \to C_{c}(\overline{N} \rtimes G)$ be the map defined by
\begin{eqnarray}
\label{subgroup}
 \varPhi(fU_{g})(x,h) &:=& \left\{\begin{array}{ll}
                        f(x)&~\text{if}~ g=h ,\\
                         0 & \text{otherwise}.
                        \end{array} \right.
\end{eqnarray}
for $f \in C_{c}(\overline{N})$ and $g \in G$. Here $\{U_{g}:g \in G\}$ denotes the canonical unitaries (corresponding to the group elements) in the multiplier algebra of $C_{0}(\overline{N}) \rtimes G$. Then $\varPhi$ extends to an isomorphism from $C_{0}(\overline{N}) \rtimes G$ onto $C^{*}(\overline{N} \rtimes G)$ (Cf. Corollary 2.3.19, Page 34, \cite{Ren_book}).

 Let $p:=1_{\overline{M}} \in C_{c}(\overline{N}) \subset C_{0}(\overline{N}) \rtimes G$ where $1_{\overline{M}}$ is the characteristic function associated to the compact open subset $\overline{M}$. Note that $\varPhi(1_{\overline{M}})=1_{\overline{M} \times \{e\}}$. 

\begin{ppsn}
\label{Corners}
The full corner $p (C_{0}(\overline{N}) \rtimes G)p$ is isomorphic to $\mathfrak{A}[N \rtimes H,M]$. Here the projection  $p$ is given by $p=1_{\overline{M}}$.
\end{ppsn}
\textit{Proof.} Let $i:C_{c}(\overline{N} \rtimes G|_{\overline{M}}) \to C_{c}(\overline{N} \rtimes G)$ be the natural inclusion. It is easy to verify that the image of $i$ is $1_{\overline{M} \times \{e\}}C_{c}(\overline{N} \rtimes G)1_{\overline{M} \times \{e\}}$. Now from Proposition \ref{Isometric embedding}, it follows that $C^{*}(\overline{N} \rtimes G|_{\overline{M}})$ is isomorphic to $1_{\overline{M} \times \{e\}}C^{*}(\overline{N} \rtimes G)1_{\overline{M} \times \{e\}}$. But we have the isomorphism $\varPhi:C_{0}(\overline{N})\rtimes G \to C^{*}(\overline{N} \rtimes G)$ with $\varPhi(1_{\overline{M}})=1_{\overline{M} \times \{e\}}$. Hence $\mathfrak{A}[N \rtimes H,M]$ is isomorphic to the corner $1_{\overline{M}}(C_{0}(\overline{N}) \rtimes G)1_{\overline{M}}$.

Let $A=C_{0}(\overline{N}) \rtimes G$. Then $ApA$ is an ideal in $A$ containing $p=1_{\overline{M}}$. Note that for every $g \in G$, $x_{g}:=U_{g}1_{\overline{M}}1_{\overline{M}} \in ApA$.  Hence $1_{g\overline{M}}= U_{g}1_{\overline{M}}U_{g}^{*}=x_{g}x_{g}^{*} \in ApA$. Hence for every $g \in G$, $1_{g. \overline{M}} \in ApA$. Thus $1_{a^{-1}\overline{M}a} \in ApA$ for every $a \in P$. Thus we have $C_{c}(\overline{N}) \subset ApA$ (See Remark \ref{compact cover}) and hence $C_{0}(\overline{N}) \subset ApA$. As a consequence we have $ApA=C_{0}(\overline{N}) \rtimes G$. Thus the projection $p$ is full. This completes the proof. \hfill $\Box$

\begin{rmrk}
\label{compact cover}
 If $K \subset \overline{N}$ is compact then there exists $b \in P$ such that $K \subset b^{-1}\overline{M}b$. For $\{a^{-1}\overline{M}a:~a \in P\}$ is an open cover of $\overline{N}$. Thus there exists $a_{1},a_{2},\cdots,a_{n} \in P$ such that $ K \subset \bigcup_{i=1}^{n}a_{i}^{-1}\overline{M}a_{i}$. Choose $ b \in \bigcap_{i=1}^{n}Pa_{i}$. Then for every $i$, $a_{i}^{-1}\overline{M}a_{i} \subset b^{-1}\overline{M}b$. (Reason: $M$ is dense in $\overline{M}$ and $ba_{i}^{-1} \in P$). Hence $K \subset b^{-1}\overline{M}b$.
\end{rmrk}

\begin{rmrk}
 Using the second half of Proposition \ref{Isometric embedding}, it can be shown that the $C^{*}$-algebra $\mathfrak{A}_{red}[N \rtimes H,M]$ is isomorphic to the full corner $1_{\overline{M}}(C_{0}(\overline{N}) \rtimes_{red} G)1_{\overline{M}}$. We leave the details to the reader.
\end{rmrk}

\section{ Cuntz-Li Duality theorem}
 The purpose of this section is to establish a duality result for the $C^{*}$-algebra associated to Examples \ref{Exel} and \ref{matrix group}. This is analogous to the duality result obtained in 
\cite{Cuntz-Li-1} for the ring $C^{*}$-algebra associated to the ring of integers in a number field. The proof is really a step by step adaptation of the arguments in \cite{Cuntz-Li-1} to our situation.

Let $\Gamma \subset GL_{n}(\mathbb{Q})$ be a subgroup and let $\Gamma_{+}:=\{\gamma \in \Gamma:\gamma \in M_{n}(\mathbb{Z})\}$. Assume that the following holds.
\begin{enumerate}
 \item The group $\Gamma=\Gamma_{+}\Gamma_{+}^{-1}=\Gamma_{+}^{-1}\Gamma_{+}$.
 \item The intersections $\displaystyle\bigcap_{\gamma \in \Gamma_{+}}\gamma \mathbb{Z}^{n}=\bigcap_{\gamma \in \Gamma_{+}} \gamma^{t}\mathbb{Z}^{n}=\{0\}$.
\end{enumerate}
Let $\Gamma^{op}:=\{\gamma^{t}:\gamma \in \Gamma\}$. Then $\Gamma^{op}$ is a subgroup of $GL_{n}(\mathbb{Q})$. Also $\Gamma$ satisfies $(1)$ and $(2)$ if and only if $\Gamma^{op}$ satisfies $(1)$ and $(2)$. If $\Gamma$  contains the non-zero scalars then $(1)$ and $(2)$ are  satisfied.

For the rest of this section, we let $\Gamma$ be a subgroup of $GL_{n}(\mathbb{Q})$ which satisfies $(1)$ and $(2)$.
The group $\Gamma$ acts on $\mathbb{Q}^{n}$ by left multiplication. Let $N_{\Gamma}:=\bigcup_{\gamma \in \Gamma_{+}}\gamma^{-1}\mathbb{Z}^{n}$. Then by Lemma \ref{The subgroup N}, it follows that $N_{\Gamma}$ is a subgroup of $\mathbb{Q}^{n}$ and $\Gamma$ leaves $N_{\Gamma}$ invariant. Consider the semidirect product $N_{\Gamma} \rtimes \Gamma$. Then the pair $(N_{\Gamma} \rtimes \Gamma,\mathbb{Z}^{n})$ satisfies the hypotheses (C1), (C2) and (C3). Let us denote the $C^{*}$-algebra $\mathfrak{A}[N_{\Gamma} \rtimes \Gamma, \mathbb{Z}^{n}]$ by $\displaystyle \mathfrak{A}_{\Gamma}$.

Note that $N_{\Gamma} \rtimes \Gamma$ acts on $\mathbb{R}^{n}$ on the right as follows. For $\xi \in \mathbb{R}^{n}$ and $(v,\gamma) \in N_{\Gamma} \rtimes \Gamma$, let $\xi.(v,\gamma)=\gamma^{-1}(\xi-v)$. This right action of $N_{\Gamma} \rtimes \Gamma$ on $\mathbb{R}^{n}$ gives rise to a left action  of $N_{\Gamma}\rtimes \Gamma$ on $C_{0}(\mathbb{R}^{n})$ as follows. For $g \in N_{\Gamma} \rtimes \Gamma$ and $f \in C_{0}(\mathbb{R}^{n})$, let $(g.f)(x)=f(x.g)$. 

 The main theorem of this section is the following.

\begin{thm}
\label{duality theorems}
 The $C^{*}$-algebras $\displaystyle \mathfrak{A}_{\Gamma^{op}}$ and $C_{0}(\mathbb{R}^{n})\rtimes (N_{\Gamma} \rtimes \Gamma)$ are Morita-equivalent.
\end{thm}

To prove this we need a bit of preparation. If $\gamma \in \Gamma_{+}$, then $\gamma$ leaves $\mathbb{Z}^{n}$ invariant and induces a map on the quotient $\displaystyle \frac{N_{\Gamma}}{\mathbb{Z}^{n}}$ which we still denote by $\gamma$. Let 

\[
 \overline{N_{\Gamma}}:=\{(z_{\gamma})_{\gamma \in \Gamma_{+}} \in \displaystyle \prod_{\gamma \in \Gamma_{+}}\frac{N_{\Gamma}}{\mathbb{Z}^{n}}:~\delta z_{\gamma \delta}=z_{\gamma} \mbox{~for every~}\gamma,\delta \in \Gamma_{+}\}
\]
We give $\displaystyle \frac{N_{\Gamma}}{\mathbb{Z}^{n}}$ the discrete topology. The abelian group $\overline{N_{\Gamma}}$  is given the subspace topology inherited from the product topology on $\displaystyle \prod_{\gamma \in \Gamma_{+}}\frac{N_{\Gamma}}{\mathbb{Z}^{n}}$. The topological group $\overline{N_{\Gamma}}$ is Hausdorff. 

Now we  describe the action of $\Gamma_{+}$ on $\overline{N_{\Gamma}}$. Let $\gamma \in \Gamma_{+}$ and $z \in \overline{N_{\Gamma}}$ be given. For $\delta \in \Gamma_{+}$, choose $\alpha, \beta \in \Gamma_{+}$ such that $\gamma \alpha = \delta\beta$.
Let $(\gamma.z)_{\delta}=\beta z_{\alpha}$. It is easily verifiable that $\gamma$ is a homeomorphism. The inverse of $\gamma$ is given by $(\gamma^{-1}z)_{\delta}=z_{\gamma\delta}$. This way $\Gamma_{+}$ acts on $\overline{N_{\Gamma}}$ and induces an action of $\Gamma$ on $\overline{N_{\Gamma}}$.
\begin{ppsn}
\label{Model}
 We have the following.
\begin{enumerate}
 \item The map $N_{\Gamma} \ni v \to (\gamma^{-1}v)_{\gamma \in \Gamma_{+}} \in \overline{N_{\Gamma}}$ is injective and is $\Gamma$-equivariant. Moreover, when $N_{\Gamma}$ is viewed as a subset of $\overline{N_{\Gamma}}$ via this embedding, $N_{\Gamma}$ is dense in $\overline{N}_{\Gamma}$.
 \item Let $\overline{M_{\Gamma}}:=\{z\in \overline{N_{\Gamma}}:z_{e}=0\}$ is a compact open subgroup of $\overline{N}_{\Gamma}$. Also  the intersection $\overline{M_{\Gamma}} \cap N_{\Gamma}=\mathbb{Z}^{n}$.
Hence $\mathbb{Z}^{n}$ is dense in $\overline{M_{\Gamma}}$.
\item  Also $\overline{N_{\Gamma}}=\bigcup_{\gamma \in \Gamma_{+}}\gamma^{-1}\overline{M_{\Gamma}}$. As a consequence, $\overline{N_{\Gamma}}$ is locally compact.
 \end{enumerate}
\end{ppsn}
\textit{Proof.} The fact that $v \to (\gamma^{-1}v)_{\gamma}$ is injective follows from the assumption that $\bigcap_{\gamma \in \Gamma_{+}}\gamma \mathbb{Z}^{n}=\{0\}$. Let $\gamma \in \Gamma_{+}$ and $v \in N_{\Gamma}$ be given. Let us denote the image of $v$ in $\overline{N_{\Gamma}}$ by $\tilde{v}$. We need to show that for $\delta \in \Gamma_{+}$, the $\delta^{\mbox{th}}$ co-ordinate of $\gamma.\tilde{v}$ is $\delta^{-1}\gamma v$. Choose $\alpha$ and $\beta$ in $\Gamma_{+}$ such that $\gamma \alpha =\delta \beta$. Then by definition $(\gamma.\tilde{v})_{\delta}=\beta \alpha^{-1}v=\delta^{-1}\gamma v$. Thus we have shown that the embedding $N_{\Gamma} \ni v \to (\gamma^{-1}v)_{\gamma \in \Gamma_{+}} \in \overline{N_{\Gamma}}$ is $\Gamma_{+}$ -equivariant and consequently is $\Gamma$ -equivariant.

For $\gamma \in \Gamma_{+}$ and $v \in N_{\Gamma}$, let \[
                                                         U_{\gamma,v}:=\{z \in \overline{N_{\Gamma}}:z_{\gamma} \equiv v \mod \mathbb{Z}^{n} \}.
                                                     \]
 Clearly the collection $\{U_{\gamma,v}:\gamma \in \Gamma_{+}, v \in N_{\Gamma} \}$ forms a basis for $\overline{N_{\Gamma}}$. Note that $\gamma.v \in U_{\gamma,v}$. Thus $N_{\Gamma}$ is dense in $\overline{N_{\Gamma}}$.

For $\gamma \in \Gamma_{+}$, let $N_{\gamma}:=\gamma^{-1}\mathbb{Z}^{n}$. Note that for $\gamma \in \Gamma_{+}$, $ \frac{N_{\gamma}}{\mathbb{Z}^{n}}$ is finite. Now observe that $\overline{M}_{\Gamma}=\overline{N_{\Gamma}} \cap  \prod_{\gamma}\frac{N_{\gamma}}{\mathbb{Z}^{n}}$. Thus $\overline{M_{\Gamma}}$ is compact. Since the projection onto the $e^{\mbox{th}}$ co-ordinate is a continuous homomorphism, it follows that $\overline{M_{\Gamma}}$ is an open subgroup. The equality $\overline{M_{\Gamma}} \cap N_{\Gamma}=\mathbb{Z}^{n}$ is obvious.

Let $z \in \overline{N_{\Gamma}}$ be given. Since $N_{\Gamma}=\bigcup_{\gamma \in \Gamma_{+}}\gamma^{-1}\mathbb{Z}^{n}$, it follows that there exists $\gamma \in \Gamma_{+}$ such that $\gamma z_{e}=0$. Then $\gamma.z \in \overline{M_{\Gamma}}$. Thus $\overline{N_{\Gamma}}=\bigcup_{\gamma \in \Gamma_{+}}\gamma^{-1}\overline{M_{\Gamma}}$. As $\overline{N_{\Gamma}}$ is a union of compact open subsets, it follows that $\overline{N_{\Gamma}}$ is locally compact. This completes the proof. \hfill $\Box$

 Let $\overline{N^{'}}$ and $\overline{M^{'}}$ be the groups considered in Section \ref{dilation} applied to the pair $(N_{\Gamma} \rtimes \Gamma, \mathbb{Z}^{n})$. Let us now convince ourselves that the pair $(\overline{N^{'}}, \overline{M^{'}})$ is  $\Gamma$-equivariantly isomorphic to the pair $(\overline{N_{\Gamma}} , \overline{M_{\Gamma}})$. Let $\gamma, \delta \in \Gamma_{+}$ be given. 

Denote the quotient map $\displaystyle \mathbb{Z}^{n} \to \frac{\mathbb{Z}^{n}}{\gamma \mathbb{Z}^{n}}$ by $q_{\gamma}$. Then $q_{\gamma}$ descends to  a map $\displaystyle \frac{\mathbb{Z}^{n}}{\gamma \delta \mathbb{Z}^{n}} \to \frac{\mathbb{Z}^{n}}{\gamma \mathbb{Z}^{n}}$ which we  denote by $q_{\gamma,\delta}$. Multiplication by $\gamma^{-1}$ maps $\mathbb{Z}^{n}$ injectively onto $\gamma^{-1}\mathbb{Z}^{n}$ and takes $\gamma \mathbb{Z}^{n}$ onto $\mathbb{Z}^{n}$. We denote the resulting isomorphism from $\displaystyle \frac{\mathbb{Z}^{n}}{\gamma \mathbb{Z}^{n}} \to \frac{\gamma^{-1}\mathbb{Z}^{n}}{\mathbb{Z}^{n}}$ again by $\gamma^{-1}$. Then we have the following commutative diagram where the vertical arrows are isomorphisms.
\begin{equation}
\label{commutative}
\begin{tikzpicture}[node distance=3cm, auto]
  \node (P) {$\displaystyle \frac{\mathbb{Z}^{n}}{\gamma \delta \mathbb{Z}^{n}}$};
  \node (B) [right of=P] {$\displaystyle \frac{\mathbb{Z}^{n}}{\gamma \mathbb{Z}^{n}}$};
  \node (A) [below of=P] {$\displaystyle \frac{(\gamma \delta)^{-1}\mathbb{Z}^{n}}{\mathbb{Z}^{n}}$};
  \node (C) [below of=B] {$\displaystyle \frac{\gamma ^{-1}\mathbb{Z}^{n}}{\mathbb{Z}^{n}}$};
   \draw[->] (P) to node {$q_{\gamma,\delta}$} (B);
  \draw[->] (P) to node [swap] {$(\gamma \delta)^{-1}$} (A);
  \draw[->] (A) to node [swap] {$\delta$} (C);
  \draw[->] (B) to node {$\gamma^{-1}$} (C);
\end{tikzpicture}
\end{equation}
Recall that 
\begin{align*}
 \overline{M^{'}}&=\{(z_{\gamma})_{\gamma \in \Gamma_{+}} \in \prod_{\gamma \in \Gamma_{+}} \frac{\mathbb{Z}^{n}}{\gamma \mathbb{Z}^{n}}: ~q_{\gamma,\delta}(z_{\gamma \delta}) = z_{\gamma} \} \\
 \overline{M_{\Gamma}}&=\{(z_{\gamma})_{\gamma \in \Gamma_{+}} \in \prod_{\gamma \in \Gamma_{+}} \frac{\gamma^{-1}\mathbb{Z}^{n}}{\mathbb{Z}^{n}}: \delta z_{\gamma \delta} = z_{\gamma} \} 
\end{align*}
Let $i:\mathbb{Z}^{n} \to \overline{M^{'}}$ be the embedding given by $i(v)=(v)_{\gamma \in \Gamma_{+}}$ and $j:\mathbb{Z}^{n} \to \overline{M_{\Gamma}}$ be the embedding described in Proposition \ref{Model}. Then $j(v)=(\gamma^{-1}v)_{\gamma \in \Gamma_{+}}$ for $v \in \mathbb{Z}^{n}$. Now the commutative diagram \ref{commutative} implies that the map $\varphi:\overline{M^{'}} \to \overline{M_{\Gamma}}$ given by $\varphi((z_{\gamma}))=(\gamma^{-1}z_{\gamma})$ is an isomorphism and $\varphi(i(v))=j(v)$ for $v \in \mathbb{Z}^{n}$. It is also clear that $\varphi$ is a homeomorphism.

\textit{Claim:} $\varphi$ is $\Gamma_{+}$-equivariant. First the embeddings $i$ and $j$ are $\Gamma_{+}$-equivariant. Since $\varphi \circ i =j$, it follows that $\varphi(\gamma.i(v))=\gamma.\varphi(i(v))$ if $\gamma \in \Gamma_{+}$ and $v \in \mathbb{Z}^{n}$. Since $i(\mathbb{Z}^{n})$ is dense in $\overline{M^{'}}$ ( and the maps involved are continuous ), it follows that $\varphi(\gamma.x)=\gamma.\varphi(x)$ for $x \in \overline{M^{'}}$ and $\gamma \in \Gamma_{+}$.

Now since $\overline{N_{\Gamma}}=\bigcup_{\gamma \in \Gamma_{+}}\gamma^{-1}\overline{M_{\Gamma}}$ and $\overline{N^{'}}=\bigcup_{\gamma \in \Gamma_{+}}\gamma^{-1}\overline{M^{'}}$, it follows from the universal property, as explained in Section 6 (item \ref{universal property}), that the map $\gamma^{-1}x \to \gamma^{-1}\varphi(x)$ (with $x \in \overline{M^{'}}$) extends to a $\Gamma$-equivariant isomorphism from $\overline{N^{'}} \to \overline{N_{\Gamma}}$. \hfill $\Box$ \newline

Now we describe the Pontryagin dual of the discrete group $N_{\Gamma}$. For $x,\xi \in \mathbb{R}^{n}$, let $<x,\xi>:=x^{t}\xi$ . If $x,\xi \in \mathbb{R}^{n}$, we let $\chi_{\xi}(x)=e^{2\pi i <x,\xi>}$. We identity $\mathbb{R}^{n}$ with $\widehat{\mathbb{R}^{n}}$ via the map $\xi \to \chi_{\xi}$. If $\xi \in \mathbb{R}^{n}$, restricting $\chi_{\xi}$ to $N_{\Gamma}$ gives a character of $N_{\Gamma}$. Moreover the map $\mathbb{R}^{n} \ni \xi \to \chi_{\xi} \in \widehat{N_{\Gamma}}$ is continuous.
 
Let $z \in \overline{N}_{\Gamma^{op}}$ be given. Let $\chi_{z}:N_{\Gamma}\to \mathbb{T}$ be defined as follows. For $x \in \gamma^{-1}\mathbb{Z}^{n}$ for some $\gamma \in \Gamma_{+}$, , let $\chi_{z}(x)=e^{2\pi i <\gamma x,z_{\gamma}>}=e^{2\pi i <x,\gamma^{t}z_{\gamma}>}$. It is easy to verify that $\chi_{z}$ is well defined and $\chi_{z}$ is a character of $N_{\Gamma}$. Clearly $\overline{N_{\Gamma^{op}}} \ni z \to \chi_{z} \in \widehat{N_{\Gamma}}$ is continuous. Note that if $z \in N_{\Gamma^{op}}$ and $x \in N_{\Gamma}$ then $\chi_{z}(x)=e^{2 \pi i <x,z>}$.

\begin{ppsn}
\label{dual}
 The map $\Psi:\mathbb{R}^{n} \times \overline{N}_{\Gamma^{op}} \to \widehat{N_{\Gamma}}$ defined by \[\Psi(\xi,z)=\chi_{\xi}\chi_{-z}\] is a surjective homomorphism with kernel $\Delta=\{(x,x):x \in N_{\Gamma^{op}}\}$. The induced map $\widetilde{\Psi}:\displaystyle \frac{\mathbb{R}^{n} \times \overline{N}_{\Gamma^{op}}}{\Delta} \to \widehat{N_{\Gamma}}$ is a topological isomorphism.
\end{ppsn}
\textit{Proof.} Clearly $\Psi$ is a continuous group homomorphism and $\Psi(\Delta)=\{1\}$. Now let us show that the kernel of $\Psi$ is $\Delta$. Let $(\xi,z)$ be such that $\Psi(\xi,z)=1$. Then for every $\gamma \in \Gamma_{+}$ and $x \in \mathbb{Z}^{n}$, we have 
\begin{align*}
 1&=\chi_{\xi}(\gamma^{-1}x)\chi_{-z}(\gamma^{-1}x) \\
  &=e^{2\pi i <x,(\gamma^{t})^{-1}\xi>}e^{-2 \pi i <x,z_{\gamma}>} \\
  &=e^{2 \pi i<x,(\gamma^{t})^{-1}\xi-z_{\gamma}>}
\end{align*}
Thus for every $\gamma \in \Gamma_{+}$, we have $z_{\gamma}-(\gamma^{t})^{-1}\xi \in \mathbb{Z}^{n}$. In other words, we have $\xi \in N_{\Gamma^{op}}$ and $z= \xi$ in $\overline{N}_{\Gamma^{op}}$. Hence $(\xi,z) \in \Delta$. Thus we have shown that the kernel of $\Psi$ is $\Delta$ which implies that $\widetilde{\Psi}$ is one-one.

Next we claim $\frac{\mathbb{R}^{n} \times \overline{N}_{\Gamma^{op}}}{\Delta}$ is compact. Let $\lambda:\mathbb{R}^{n} \times \overline{N}_{\Gamma^{op}} \to \frac{\mathbb{R}^{n} \times \overline{N}_{\Gamma^{op}}}{\Delta}$ be the quotient map. We also write $\lambda(\xi,z)$ as $[(\xi,z)]$. We claim that $\lambda([0,1]^{n} \times \overline{M}_{\Gamma^{op}})=\frac{\mathbb{R}^{n} \times \overline{N}_{\Gamma^{op}}}{\Delta}$. This will prove that $\frac{\mathbb{R}^{n} \times \overline{N}_{\Gamma^{op}}}{\Delta}$ is compact. 

Let $[(\xi,z)]$ be an element in the quotient $\frac{\mathbb{R}^{n} \times \overline{N}_{\Gamma^{op}}}{\Delta}$. Choose $v \in \mathbb{Z}^{n}$ and $\gamma \in \Gamma_{+}$ such that $z_{e}\equiv (\gamma^{t})^{-1}v$. Then $[(\xi,z)]=[(\xi-(\gamma^{t})^{-1}v,z-(\gamma^{t})^{-1}v)]$. Choose $w \in \mathbb{Z}^{n}$ such that $\xi-(\gamma^{t})^{-1}v-w \in [0,1]^{n}$. Let $\xi^{'}=\xi-(\gamma^{t})^{-1}v-w$ and $z^{'}=z-(\gamma^{t})^{-1}v-w $.
 Then $\xi^{'} \in [0,1]^{n}$ and $z^{'} \in \overline{M}_{\Gamma^{op}}$. Moreover $\lambda(\xi,z)=\lambda(\xi^{'},z^{'})$. 
Thus the image of $[0,1]^{n} \times \overline{M}_{\Gamma^{op}}$ under $\lambda$ is $\frac{\mathbb{R}^{n} \times \overline{N}_{\Gamma^{op}}}{\Delta}$.

The image of $\widetilde{\Psi}$ is a compact subgroup of $\widehat{N_{\Gamma}}$ and it separates points of $N_{\Gamma}$ ( The image of $\mathbb{R}^{n} \times \{0\}$ under $\Psi$ separates points of $N_{\Gamma}$). Hence $\widetilde{\Psi}$ is onto.
Since $\frac{\mathbb{R}^{n} \times \overline{N}_{\Gamma^{op}}}{\Delta}$ is compact, it follows that $\widetilde{\Psi}$ is a topological isomorphism. This completes the proof. \hfill $\Box$
                                                                                                                                                                                                                                                                                                                                                                                                                                                                                                                                                                                                          
Consider the semidirect product $\mathbb{R}^{n} \rtimes \Gamma^{op}$ where $\Gamma^{op}$ acts on $\mathbb{R}^{n}$ by left multiplication. The semidirect product $\mathbb{R}^{n}\rtimes \Gamma^{op}$ acts on 
$\widehat{N_{\Gamma}}=\frac{\mathbb{R}^{n}\times \overline{N_{\Gamma^{op}}}}{\Delta}$ on the right as follows. For $[(\xi,z)] \in \widehat{N_{\Gamma}}$ and $(v,\gamma) \in \mathbb{R}^{n}\rtimes \Gamma^{op}$, let $[(\xi,z)].(v,\gamma)=[(\gamma^{-1}(\xi+v),\gamma^{-1}z)]$. This right action of $\mathbb{R}^{n} \rtimes \Gamma^{op}$ on $\widehat{N_{\Gamma}}$ induces a left action of $\mathbb{R}^{n} \rtimes \Gamma^{op}$ on $C^{*}(N_{\Gamma}) \cong C(\widehat{N_{\Gamma}})$. 

 The crossed product $C^{*}(N_{\Gamma}) \rtimes (\mathbb{R}^{n} \rtimes \Gamma^{op})$ is isomorphic to the iterated crossed product 
$(C^{*}(N_{\Gamma}) \rtimes \mathbb{R}^{n}) \rtimes \Gamma^{op}$. (Cf. Proposition 3.11, Page 87, \cite{Dana-Williams}.) But then the map $\Gamma \ni \gamma \to (\gamma^{t})^{-1} \in \Gamma^{op}$ is an isomorphism. Thus the crossed product $(C^{*}(N_{\Gamma}) \rtimes \mathbb{R}^{n}) \rtimes \Gamma^{op} \cong (C^{*}(N_{\Gamma}) \rtimes \mathbb{R}^{n}) \rtimes \Gamma$.

Let us fix notations. Let $\tau$ be the action of $\mathbb{R}^{n}$ on $C^{*}(N_{\Gamma})$. Let $\beta$ be the action of $\Gamma$ on $C^{*}(N_{\Gamma}) \cong C(\widehat{N_{\Gamma}})$, induced by the action of $\Gamma^{op}$ and the identification $\Gamma \cong \Gamma^{op}$. For $v \in N_{\Gamma}$, $\xi \in \mathbb{R}^{n}$ and $\gamma \in \Gamma$, it is easy to verify the following.
\begin{align*}
 \tau_{\xi}(\delta_{v})&=e^{-2 \pi i <\xi,v>}\delta_{v}, \\
 \beta_{\gamma}(\delta_{v})&=\delta_{\gamma v}
\end{align*}
where $\{\delta_{v}:v \in N_{\Gamma}\}$ denotes the canonical unitaries of $C^{*}(N_{\Gamma})$. The action of $\Gamma^{op}$ on $C^{*}(N_{\Gamma}) \rtimes \mathbb{R}^{n}$, induces an action of $\Gamma$  ( via the identification $\Gamma \ni \gamma \to (\gamma^{t})^{-1}$) and let us denote it by $\widetilde{\beta}$. For $\gamma \in \Gamma$, and $f \in C_{c}(\mathbb{R}^{n},C^{*}(N_{\Gamma}))$, we have 
\[
 \widetilde{\beta_{\gamma}}(f)(x)= |\det(\gamma)|\beta_{\gamma}(f(\gamma^{t}x)).
\]

Now consider the crossed product $C_{0}(\mathbb{R}^{n}) \rtimes (N_{\Gamma} \rtimes \Gamma) \cong C^{*}(\mathbb{R}^{n}) \rtimes (N_{\Gamma} \rtimes \Gamma)$. Let us denote the action of $N_{\Gamma}$ and $\Gamma$ on $C^{*}(\mathbb{R}^{n})$ by $\sigma$ and $\alpha$. For $v \in N_{\Gamma}$, $\gamma \in \Gamma$ and $f \in C_{c}(\mathbb{R}^{n})$, we have 
\begin{align*}
 (\sigma_{v}f)(\xi)&=e^{2 \pi i <\xi,v>}f(\xi), \\
 (\alpha_{\gamma}f)(\xi)&=|\det(\gamma)|f(\gamma^{t}\xi).
\end{align*}
Denote the action of $\Gamma$ on $C^{*}(\mathbb{R}^{n}) \rtimes N_{\Gamma}$ by $\widetilde{\alpha}$. For $\gamma \in \Gamma$, $v \in N_{\Gamma}$ and $f \in C^{*}(\mathbb{R}^{n})$, one has \[
                                                                                \widetilde{\alpha_{\gamma}}(f\delta_{v})=\alpha_{\gamma}(f)\delta_{\gamma v}.
                                                                                                                                                                                            \]

Let us recall the following lemma which is   Lemma 4.3 in \cite{Cuntz-Li-1}.

\begin{lmma}[\cite{Cuntz-Li-1}]
\label{duality lemma}
Let $G$ be a locally compact abelian group and $H$ be a subgroup of the Pontryagin dual $\widehat{G}$. Endow $H$ with the discrete topology. Let $\sigma$ be the action of $H$ on $C^{*}(G)$ and $\tau$ be the action of $G$ on $C^{*}(H)$ given by $\sigma_{h}(f)=[g \to h(g)f(g)]$ and $\tau_{g}(\tilde{f})=[h \to h(-g)\tilde{f}(h)]$. Then the map $\phi:C_{c}(H,C_{c}(G)) \to C_{c}(G,C_{c}(H))$ defined by $\phi(f)(g)(h)=h(-g)f(h)(g)$ extends to an isomorphism between $C^{*}(G) \rtimes_{\sigma} H$ and $C^{*}(H) \rtimes_{\tau} G$.
 \end{lmma}

We are now ready to prove the following proposition. 

\begin{ppsn}
\label{duality}
 The crossed products $C_{0}(\mathbb{R}^{n})\rtimes (N_{\Gamma} \rtimes \Gamma)$ and $C(\widehat{N_{\Gamma}})\rtimes (\mathbb{R}^{n} \rtimes \Gamma^{op})$ are isomorphic.
\end{ppsn}
\textit{Proof.} It is enough to show that the crossed products $(C^{*}(\mathbb{R}^{n})\rtimes_{\sigma} N_{\Gamma}) \rtimes_{\widetilde{\alpha}} \Gamma$ and  $(C^{*}(N_{\Gamma}) \rtimes_{\tau} \mathbb{R}^{n}) \rtimes_{\widetilde{\beta}} \Gamma$ are isomorphic. We show that $C^{*}(\mathbb{R}^{n}) \rtimes_{\sigma} N_{\Gamma}$ and $C^{*}(N_{\Gamma}) \rtimes_{\tau} \mathbb{R}^{n}$ are $\Gamma$-equivariantly isomorphic. Then the isomorphism between the crossed products will follow.

Identify $\mathbb{R}^{n}$ with $\widehat{\mathbb{R}^{n}}$ via the map $\xi \to \chi_{\xi}$. (Recall that $\chi_{\xi}$ is the character given by $\chi_{\xi}(x)=e^{2 \pi i <x,\xi>}$. )  Consider $N_{\Gamma}$ as a subgroup of $\widehat{\mathbb{R}^{n}}$ via the natural inclusion $N_{\Gamma} \subset \mathbb{R}^{n}$. Note that the action $\sigma$ of $N_{\Gamma}$ on $C^{*}(\mathbb{R}^{n})$ and $\tau$ of $\mathbb{R}^{n}$ on $C^{*}(N_{\Gamma})$ are exactly as in Lemma \ref{duality lemma}. 

Thus Lemma \ref{duality lemma} implies that $C^{*}(\mathbb{R}^{n}) \rtimes_{\sigma} N_{\Gamma} \cong C^{*}(N_{\Gamma}) \rtimes_{\tau} \mathbb{R}^{n}$. Let $\phi:C^{*}(\mathbb{R}^{n}) \rtimes_{\sigma} N_{\Gamma} \to C^{*}(N_{\Gamma}) \rtimes_{\tau} \mathbb{R}^{n}$ be the isomorphism prescribed by Lemma \ref{duality lemma}. We claim $\phi$ is $\Gamma$-equivariant.  First note that $\phi(f\delta_{v})(\xi)=e^{-2 \pi i <\xi,v>}f(\xi)\delta_{v}$ for $f \in C_{c}(\mathbb{R}^{n})$ and $v \in N_{\Gamma}$.

Let $\gamma \in \Gamma$ be given. Now observe that \begin{align*}
\widetilde{\beta_{\gamma}}(\phi(f\delta_{v}))(\xi)&=|\det(\gamma)|\beta_{\gamma}(\phi(f\delta_{v})(\gamma^{t}\xi)) \\
                                                &=|\det(\gamma)|e^{-2 \pi i <\gamma^{t}\xi,v>}f(\gamma^{t}\xi)\delta_{\gamma v}\\
                                                &=|\det(\gamma)|e^{-2 \pi i <\xi,\gamma v>}f(\gamma^{t}\xi)\delta_{\gamma v}.
                                                                                       \end{align*}
On the other hand, observe that \begin{align*}
                                 \phi(\widetilde{\alpha_{\gamma}}(f\delta_{v}))(\xi)&=\phi(\alpha_{\gamma}(f)\delta_{\gamma v})(\xi) \\
                                &= e^{-2 \pi i <\xi , \gamma v>}\alpha_{\gamma}(f)(\xi)\delta_{\gamma v} \\
                                &=e^{-2 \pi i <\xi,\gamma v>}|\det(\gamma)|f(\gamma^{t}\xi)\delta_{\gamma v}
                                \end{align*}

Hence for every $\gamma \in \Gamma$, $\widetilde{\beta_{\gamma}}\phi(f\delta_{v})=\phi \widetilde{\alpha_{\gamma}}(f \delta_{v})$. Since  $\{f \delta_{v}: f \in C_{c}(\mathbb{R}^{n}) , v \in N_{\Gamma} \}$ is total in $C^{*}(\mathbb{R}^{n}) \rtimes_{\sigma} N_{\Gamma}$, it follows that for every $\gamma$, $\widetilde{\beta_{\gamma}}\phi=\phi \widetilde{\alpha_{\gamma}}$. In other words, $\phi$ is $\Gamma$-equivariant.
This completes the proof. \hfill $\Box$

\textbf{Proof of Theorem \ref{duality theorems}}.
 By Corollary \ref{main corollary}, it follows that $\mathfrak{A}_{\Gamma^{op}}$ is isomorphic to the $C^{*}$-algebra of the groupoid $\widetilde{\mathcal{G}}:=\overline{N}_{\Gamma^{op}}\rtimes (N_{\Gamma^{op}} \rtimes \Gamma^{op})|_{\overline{M}_{\Gamma^{op}}}$. By Proposition \ref{duality}, it follows that $C_{0}(\mathbb{R}^{n})\rtimes (N_{\Gamma}\rtimes \Gamma)$ is isomorphic to the $C^{*}$-algebra of the groupoid $\mathcal{G}:=\widehat{N_{\Gamma}}\rtimes (\mathbb{R}^{n} \rtimes \Gamma^{op})$. We will show that $\mathcal{G}$ and $\widetilde{\mathcal{G}}$ are equivalent in the sense of \cite{MRW}.

By Proposition \ref{dual}, $\widehat{N}_{\Gamma}=\frac{\mathbb{R}^{n}\times \overline{N}_{\Gamma^{op}}}{\Delta}$ where $\Delta:=\{(x,x):x \in N_{\Gamma^{op}}\}$. Denote the quotient map $\mathbb{R}^{n}\times \overline{N}_{\Gamma^{op}} \to \frac{\mathbb{R}^{n}\times \overline{N}_{\Gamma^{op}}}{\Delta}$ by $\lambda$.   Let $X:=\lambda(\{0\} \times \overline{M}_{\Gamma^{op}})$.  Then $X$ is a closed  subset of $\mathcal{G}^{0}$ and it is easy to verify that $X$ meets each orbit of $\mathcal{G}^{0}$. Let 
\begin{displaymath}
 \mathcal{G}_{X}:=\{ \alpha \in \mathcal{G}:s(\alpha) \in X\}=s^{-1}(X)
\end{displaymath}
We claim that the (restricted) source map $s:\mathcal{G}_{X} \to X$ and the range map $r:\mathcal{G}_{X} \to \mathcal{G}^{0}$ are open. Let $U\subset \mathcal{G}$ be an open subset. Then $s(U\cap \mathcal{G}_{X})= s(U)\cap X$. Since $s:\mathcal{G}\to \mathcal{G}^{0}$ is open, it follows that $s:\mathcal{G}_{X} \to X$ is open.

Now we prove that $r:\mathcal{G}_{X} \to \mathcal{G}^{0}$ is open. It is enough to show that $r((U \times V \times \{\gamma \}) \cap \mathcal{G}_{X})$ is open whenever $U \subset \frac{\mathbb{R}^{n}\times \overline{N}_{\Gamma^{op}}}{\Delta}$ and $V \subset \mathbb{R}^{n}$ are open and $\gamma \in \Gamma^{op}$. We claim that 
\[
 r((U \times V \times \{\gamma \})\cap \mathcal{G}_{X})=U \cap \lambda(-V \times \gamma \overline{M}_{\Gamma^{op}})
\]
Let $[(\xi,z)] \in r((U \times V \times \{\gamma\})\cap \mathcal{G}_{X})$. Then there exists $([(\eta,y)],v,\gamma) \in U \times V \times \{\gamma\}$ such that $[(\eta,y)].(v,\gamma) \in X$ and $[(\xi,z)]=[(\eta,y)]$. Thus there exists $u \in N_{\Gamma^{op}}$ such that $\gamma^{-1}(\xi+v)=u$ and $\gamma^{-1}z-u=x$ for some $x \in \overline{M}_{\Gamma^{op}}$. Hence $[(\xi,z)]=[(-v,\gamma x)]$. Clearly $[(\xi,z)] \in U$. Hence $[(\xi,z)] \in U \cap \lambda(-V \times \gamma\overline{M}_{\Gamma^{op}}).$ Thus we have shown that \[r((U \times V \times \{\gamma\})\cap \mathcal{G}_{X})\subset U \cap \lambda(-V \times \gamma\overline{M}_{\Gamma^{op}}).\]

Now let $[(\xi,z)] \in U \cap \lambda(-V \times \gamma\overline{M}_{\Gamma^{op}})$. Then there exists $(v,x) \in V \times \overline{M}_{\Gamma^{op}}$ such that $[(\xi,z)]=[(-v,\gamma x)]$. This is equivalent to saying that $[(\xi,z)].(v,\gamma) \in X$. Thus $([(\xi,z)],v,\gamma) \in (U \times V \times \{\gamma\})\cap \mathcal{G}_{X}$ and $r([(\xi,z)],v,\gamma)=(\xi,z)]$. This proves that $U \cap \lambda(-V \times \gamma\overline{M}_{\Gamma^{op}}) \subset r((U \times V \times \{\gamma\})\cap \mathcal{G}_{X})$.

This proves the claim that $r((U \times V \times \{\gamma\})\cap \mathcal{G}_{X})=U \cap \lambda(-V \times \gamma \overline{M}_{\Gamma^{op}})$. Now since $\lambda$ is open and $\overline{M}_{\Gamma^{op}}$ is open, it follows that $r((U \times V \times \{\gamma \})\cap \mathcal{G}_{X})$ is open. Thus we have shown that $r:\mathcal{G}_{X} \to \mathcal{G}^{0}$ is open.

Now by Example 2.7 of \cite{MRW}, it follows that  $\mathcal{G}$ and $\mathcal{G}_{X}^{X}:=\{\alpha \in \mathcal{G}_{X}:r(\alpha) \in X\}$ are equivalent. Recall that $\widetilde{\mathcal{G}}=\overline{N}_{\Gamma^{op}} \rtimes (N_{\Gamma^{op}}\rtimes \Gamma^{op})|_{\overline{M}_{\Gamma^{op}}}$  The right action of $N_{\Gamma^{op}}\rtimes \Gamma^{op}$ on $\overline{N}_{\Gamma^{op}}$ is given by 
$x.(v,\gamma)=\gamma^{-1}(x-v)$. Let $\Phi:\widetilde{\mathcal{G}} \to \mathcal{G}_{X}^{X}$ be defined by $\Phi(x,v,\gamma)=([(0,x)],v,\gamma)$. It is easy to check that $\Phi$ is a groupoid isomorphism and it is continuous. Now we prove that $\Phi$ is a topological isomorphism.

Let $(x_{n},v_{n},\gamma)$ be a sequence in $\widetilde{\mathcal{G}}$ such that $\Phi(x_{n},v_{n},\gamma)$ converges to $([(0,x)],v,\gamma))$. First note that $x \to [(0,x)]$ is a topological embedding of $\overline{M}_{\Gamma^{op}}$ into $\widehat{N_{\Gamma}}$. Thus, it follows that $x_{n}$ converges to $x$ in $\overline{M}_{\Gamma^{op}}$. Now $\Phi(x_{n},v_{n},\gamma)$ converges to $[(0,x)],v,\gamma)$ implies that $v_{n}$ tends to $v$ in $\mathbb{R}^{n}$ and $\gamma^{-1}(x-v_{n})$ tends to $\gamma^{-1}(x-v)$ in $\overline{M}_{\Gamma^{op}}$. Hence $v_{n}$ converges to $v$ in $\overline{N}_{\Gamma^{op}}$. Thus $(v_{n},v_{n}) \to (v,v)$ in $\mathbb{R}^{n}\times \overline{N}_{\Gamma^{op}}$. But $\Delta$ is a discrete subgroup of $\mathbb{R}^{n}\times \overline{N}_{\Gamma^{op}}$. Hence $v_{n}=v$ eventually. Therefore, $(x_{n},v_{n},\gamma)\to (x,v,\gamma)$ in $\widetilde{\mathcal{G}}$.
So, $\Phi$ is a topological isomorphism. 

Since $\mathcal{G}$ and $\widetilde{\mathcal{G}}$ are equivalent in the sense of \cite{MRW}, it follows from Theorem 2.8 in \cite{MRW} that $C^{*}(\mathcal{G})$ and $C^{*}(\widetilde{\mathcal{G}})$ are Morita-equivalent. This completes the proof. \hfill $\Box$

\subsection{Examples }
We end this article by considering two examples. 

\textbf{Example 1:} First we show that the duality result for the ring $C^{*}$-algebra  associated to number fields obtained in \cite{Cuntz-Li-1} can be derived from Theorem \ref{duality theorems}.

Consider  a number field $K$ of degree $n$. Denote the ring of integers in $K$ by $O_{K}$. Let $\{w_{1},w_{2},\cdots,w_{n}\}$ be a $\mathbb{Z}$-basis for $O_{K}$. Then $\{w_{1},w_{2},\cdots,w_{n}\}$ is a $\mathbb{Q}$-basis for $K$. Identify $K$ with $\mathbb{Q}^{n}$ via the map $\beta:\mathbb{Q}^{n} \ni (x_{1},x_{2},\cdots,x_{n})^{t}\to \sum_{i=1}^{n}x_{i}w_{i} \in K$. By definition, $\beta(\mathbb{Z}^{n})=O_{K}$.

If $a \in K$, then $a$ acts on $K$ by left multiplication and is $\mathbb{Q}$-linear. Thus $a$ gives rise to a matrix with respect to the basis $\{w_{1},w_{2},\cdots,w_{n}\}$ which we denote by $\alpha(a)$. Explicitly, for $1 \leq j \leq n$, let
 \begin{equation}
\label{matrix representation}                                                                                                                                                                                                                                                      
aw_{j}:=\sum_{i=1}^{n}\alpha_{ij}(a)w_{i}.                                                                                                                                                                                                                                                      \end{equation}

 Let $\alpha(a):=(\alpha_{ij}(a))$. Then $\alpha:K \to M_{n}(\mathbb{Q})$ is an injective ring homomorphism. We also have the following equivariance. For $a \in K$ and $x \in \mathbb{Q}^{n}$,  $\beta(\alpha(a)x)=a\beta(x)$.

Let $\Gamma:=\alpha(K^{\times})$. Then $\Gamma$ is a subgroup of $GL_{n}(\mathbb{Q})$. Now the pair $(K \rtimes K^{\times},O_{K})$ is isomorphic to $(\mathbb{Q}^{n}\rtimes \Gamma, \mathbb{Z}^{n})$. Thus the ring $C^{*}$-algebra associated to $O_{K}$ is nothing but $\mathfrak{A}[\mathbb{Q}^{n}\rtimes \Gamma, \mathbb{Z}^{n}]$. Hence Theorem \ref{duality theorems} applies. The only thing that one needs to verify is  $\bigcap_{a \in O_{K}}\alpha(a)^{t}\mathbb{Z}^{n}$ is trivial. Since $\bigcap_{a \in O_{K}}aO_{K}=\{0\}$, it follows that $\bigcap_{a \in O_{K}}\alpha(a)\mathbb{Z}^{n}=\{0\}$. We produce a matrix $X$ with rational entries whose determinant is non-zero and $X\alpha(a)X^{-1}=\alpha(a)^{t}$ for every $a \in O_{K}$. Then it will follow that $\displaystyle \bigcap_{a \in O_{K}}\alpha(a)^{t}\mathbb{Z}^{n}=\{0\}$. (See also Lemma \ref{Conjugation to transpose}.)

Let $Tr:M_{n}(\mathbb{Q}) \to \mathbb{Q}$ be the usual trace and let $tr:=Tr\circ \alpha$. Denote the $n\times n$ matrix whose $(i,j)^{\mbox{th}}$ entry is $tr(w_{i}w_{j})$ by $X$. Then $X$ has determinant non-zero and its determinant is called the discriminant of the number field $K$.

\begin{lmma}
 \label{Discriminant matrix}
For every $a \in K$, $X\alpha(a)X^{-1}=\alpha(a)^{t}$.
\end{lmma}
\textit{Proof.} Fix $a \in K$. Let $Y=(tr(aw_{i}w_{j}))$. Multiplying Equation \ref{matrix representation} by $w_{k}$ and taking  trace, we get 
\begin{displaymath}
 Y_{jk}=\sum_{i=1}^{n}\alpha_{ij}(a)X_{ik}
\end{displaymath}
In other words, we have $Y=\alpha(a)^{t}X$. But $Y$ and $X$ are symmetric. Thus taking transpose, we get $Y=X \alpha(a)$. Hence $X\alpha(a)=\alpha(a)^{t}X$. This completes the proof. \hfill $\Box$

Let $\mathbb{A}_{\infty}$ denote the ring of infinite adeles associated to $K$.
\begin{thm}[\cite{Cuntz-Li-1}]
 For a number field $K$, the ring $C^{*}$-algebra $\mathfrak{A}[K \rtimes K^{\times},O_{K}]$ is Morita-equivalent to $C_{0}(\mathbb{A}_{\infty})\rtimes (K \rtimes K^{\times})$.
\end{thm}
\textit{Proof.} Note that for $\Gamma=\alpha(K^{\times})$, $N_{\Gamma}=\mathbb{Q}^{n}$ and $N_{\Gamma^{op}}=\mathbb{Q}^{n}$ ( since $\Gamma$ contains the diagonal matrices with rational entries). Thus Lemma \ref{Discriminant matrix} implies that the matrix $X=(tr(w_{i}w_{j}))$ implements an isomorphism between the dynamical systems $(\mathbb{R}^{n},\mathbb{Q}^{n} \rtimes \Gamma)$ and $(\mathbb{R}^{n}, \mathbb{Q}^{n} \rtimes \Gamma^{op})$. The map \[(\mathbb{R}^{n},\mathbb{Q}^{n} \rtimes \Gamma) \ni (\xi,(v,\gamma)) \to (X\xi,(Xv,\gamma^{t})) \in (\mathbb{R}^{n}, \mathbb{Q}^{n} \rtimes \Gamma^{op})\] is the required isomorphism. (Note that $\Gamma$ is commutative.)

 Consider the map  $\delta:\mathbb{R}^{n} \ni (x_{1},x_{2},\cdots,x_{n}) \to \sum_{i=1}x_{i}w_{i} \in \mathbb{A}_{\infty}$. Then from standard number theoretic arguments, ( for example, using Theorem 13.5 (page 70)  and Theorem 4.4 (page 110) in \cite{Gerald}), it follows that $\delta$ (together with identifications $\alpha$ and $\beta$) implements an isomorphism between  $(\mathbb{A}_{\infty},K \times K^{\rtimes})$ and $(\mathbb{R}^{n},\mathbb{Q}^{n} \rtimes \Gamma)$. Now Theorem \ref{duality theorems} yields the required result. This completes the proof. \hfill $\Box$ 

\textbf{Example 2:} Let $A$ be an $n \times n$ matrix with integer entries such that $\det(A) \neq 0$ and $\bigcap_{r=0}^{\infty}A^{r}\mathbb{Z}^{n}=\{0\}$. Let $\Gamma:=\{A^{r}: r \in \mathbb{Z} \} \cong \mathbb{Z}$. Denote the subgroup $N_{\Gamma}$ by $N_{A}$ and the Cuntz-Li algebra $\mathfrak{A}[N_{\Gamma} \rtimes \Gamma, \mathbb{Z}^{n}]$ by $\mathfrak{A}_{A}$.  Denote the transpose $A^{t}$ by $B$. Then $\Gamma^{op}=\{B^{r}:r \in \mathbb{Z}\} \cong \mathbb{Z}$. 

 We claim that the duality result is applicable to this example. The only thing that needs verification is $\bigcap_{r=0}^{\infty}B^{r}\mathbb{Z}^{n}=\{0\}$. This follows from the following lemma.

\begin{lmma}
\label{Conjugation to transpose}
 Let $A$ be a $n\times n$ matrix with integer entries and denote $A^{t}$ by $B$. Then
$\bigcap_{r=0}^{\infty}A^{r}\mathbb{Z}^{n}=\{0\}$ if and only if $\bigcap_{r=0}^{\infty}B^{r}\mathbb{Z}^{n}=\{0\}$.
\end{lmma}
\textit{Proof.} Since $A$ and $B$ are similar over $\mathbb{Q}$, it follows that there exists $Y \in GL_{n}(\mathbb{Q})$ such that $YAY^{-1}=B$. Choose a non-zero integer $m$ such that $X=mY \in M_{n}(\mathbb{Z})$. One has $XA=BX$. By induction, it follows that $XA^{r}=B^{r}X$ for every $r \geq 0$. First note that it is enough to show that $\bigcap_{r=0}^{\infty}A^{r}\mathbb{Z}^{n} \neq \{0\}$ implies $\bigcap_{r=0}^{\infty}B^{r}\mathbb{Z}^{n} \neq \{0\}$.

Suppose $v$ is a non-zero element in $\bigcap_{r=0}^{\infty}A^{r}\mathbb{Z}^{n}$. Then 
\begin{align*}
Xv &\in \bigcap_{r=0}^{\infty}XA^{r}\mathbb{Z}^{n}\\ & = \bigcap_{r=0}^{\infty}B^{r}X\mathbb{Z}^{n} \subset \bigcap_{r=0}^{\infty}B^{r}\mathbb{Z}^{n}.
\end{align*}
Since $X$ is invertible over $\mathbb{Q}$, it follows that $Xv$ is a non-zero element in $\bigcap_{r=0}^{\infty}B^{r}\mathbb{Z}^{n}$. Thus if $\bigcap_{r=0}^{\infty}A^{r}\mathbb{Z}^{n} \neq \{0\}$ then $\bigcap_{r=0}^{\infty}B^{r}\mathbb{Z}^{n} \neq \{0\}$. This completes the proof. \hfill $\Box$

Now Theorem \ref{duality theorems} and Proposition \ref{duality} implies the following proposition. 

\begin{ppsn}
\label{Morita for integer dilation}
 The $C^{*}$-algebra $\mathfrak{A}_{A^{t}}$ is Morita-equivalent to $C_{0}(\mathbb{R}^{n}) \rtimes (N_{A} \rtimes \mathbb{Z})$. Also $\mathfrak{A}_{A^{t}}$ is Morita-equivalent to $ (C^{*}(N_{A}) \rtimes \mathbb{R}^{n}) \rtimes \mathbb{Z}$. 
\end{ppsn}

Proposition \ref{Morita for integer dilation} for the case when  $n=1$ and $A=(2)$ was proved in \cite{Larsen-Li-1}. In this case, the $C^{*}$-algebra $\mathfrak{A}_{A^{t}}=\mathfrak{A}_{A}$ is the $C^{*}$-algebra $\mathcal{Q}_{2}$ considered in \cite{Larsen-Li-1}. The subgroup $\bigcup_{r=0}2^{-r}\mathbb{Z}$ is denoted $\mathbb{Z}[\frac{1}{2}]$ in \cite{Larsen-Li-1}. The Morita equivalence between $\mathcal{Q}_{2}$ and $C_{0}(\mathbb{R}) \rtimes (\mathbb{Z}[\frac{1}{2}] \rtimes (2))$ is called the $2$-adic duality theorem in \cite{Larsen-Li-1}. (Cf. Corollary 5.5 and Theorem 7.5 in \cite{Larsen-Li-1}.)

\nocite{Muhly}

\bibliography{references}
\bibliographystyle{amsalpha}
\end{document}